\documentclass[review,onefignum,onetabnum]{siamart171218}
\usepackage{array} 
\usepackage[caption=false]{subfig}
\usepackage{siunitx}
\allowdisplaybreaks
\nolinenumbers
\usepackage{tikz}
\usetikzlibrary{positioning}
\usetikzlibrary{calc}
\usetikzlibrary{chains}
\usetikzlibrary{shapes}

\usepackage{lipsum}
\usepackage{amsfonts}
\usepackage{graphicx}
\usepackage{epstopdf}
\usepackage{algorithmic}
\ifpdf
  \DeclareGraphicsExtensions{.eps,.pdf,.png,.jpg}
\else
  \DeclareGraphicsExtensions{.eps}
\fi

\newsiamremark{remark}{Remark}
\newsiamremark{hypothesis}{Hypothesis}
\crefname{hypothesis}{Hypothesis}{Hypotheses}
\newsiamthm{claim}{Claim}

\headers{Integral Equation Methods for Amphiphilic Interaction}{S.-P. P. Fu, R. J. Ryham, A. Kl{\"o}ckner, M. Wala, S. Jiang, and Y.-N. Young}

\title{Simulation of Multiscale Hydrophobic Lipid Dynamics via Efficient Integral Equation Methods\thanks{
Submitted to the editors DATE.
}}

\author{Szu-Pei P. Fu\thanks{Department of Mathematics, Fordham University, Bronx, NY 10458
  (\email{sfu17@fordham.edu}, \email{rryham@fordham.edu})}
\and Rolf J. Ryham\footnotemark[2]
\and  Andreas Kl{\"o}ckner\thanks{Department of Computer Science, University of Illinois at Urbana-Champaign, Urbana, IL 61801
  (\email{andreask@illinois.edu}, \email{wala1@illinois.edu}).}
  \and Matt Wala\footnotemark[3]
  \and Shidong~Jiang\thanks{Department of Mathematical Sciences, New Jersey Institute of Technology, Newark, NJ
  07102
  (\email{shidong.jiang@njit.edu }, \email{yyoung@njit.edu}).}
    \and Y.-N. Young\footnotemark[4]}

\usepackage{amsopn}

\ifpdf
\hypersetup{
  pdftitle={Simulation of Multiscale Hydrophobic Lipid Dynamics via Efficient Integral Equation Methods},
  pdfauthor={Szu-Pei P. Fu, Rolf J. Ryham, Andreas Kl{\"o}ckner, Matt Wala, Shidong Jiang and Yuan-Nan Young}
}
\fi

\newcommand*{\thead}[1]{\multicolumn{1}{|c|}{\bfseries #1}}

\newcommand{\KB}{k_{\mathrm{B}}}

\newcommand{\KTH}{k_{\theta}}
\newcommand{\KA}{k_{\mathrm{A}}}
\newcommand{\KBT}{\mathrm{k_BT}}

\DeclareMathOperator{\Div}{Div}
\DeclareMathOperator{\Curl}{Curl}

\everymath{\displaystyle}

\begin{document}
\maketitle

\begin{abstract}
In this paper a mathematical model for long-range, hydrophobic attraction between amphiphilic particles is 
developed to quantify the macroscopic assembly and mechanics of a lipid bilayer membrane in solvents. 
The non-local interactions between amphiphilic particles are obtained from the first domain variation of a 
hydrophobicity functional, giving rise to forces and torques (between particles) that dictate the motion of both particles
and the surrounding solvent. 
The functional minimizer (that accounts for hydrophobicity at molecular-aqueous interfaces) is a solution to a boundary value problem of the screened Laplace equation. 
We reformulate the boundary value problem as a second-kind
integral equation (SKIE), discretize the SKIE using a Nystr\"om discretization and `Quadrature by Expansion' (QBX)
and solve the resulting linear system iteratively using GMRES.
We evaluate the required layer potentials using the `GIGAQBX' fast algorithm, a variant of the
Fast Multipole Method (FMM), yielding the required particle interactions with asymptotically optimal cost.
Solving a mobility problem in Stokes flow is incorporated to obtain corresponding rigid body motion.
The simulated fluid-particle systems exhibit a variety of multiscale behaviors over both time and length:
Over short time scales, the numerical results show self-assembly for model lipid particles. 
For large system simulations, the particles form realistic configurations like micelles and bilayers. 
Over long time scales, the bilayer shapes emerging from the simulation appear to minimize a form of bending energy.
\end{abstract}

\begin{keywords}
Energy Variation, Integral Equation Method, Lipid Dynamics
\end{keywords}

\begin{AMS}
  31A10, 
  35A15, 
  92C05 
\end{AMS}

\section{Introduction}
\label{sec:intro}
In recent years, researchers have developed various macroscopic continuum formulations and
a number of numerical methods for calculating
energy minimizing and time-dependent shapes of lipid bilayer membranes, vesicles and red blood cells. While the Helfrich free energy of a lipid bilayer membrane assumes an infinitely thin membrane thickness~\cite{FriedSeguin14,Helfrich73},
many other continuum models incorporate more lipid physics \cite{Bartels,Burger13,OhtaKawasaki} and 
membrane structures \cite{DaiPromislow2013_RSPa,DaiPromislow2015_SIAMJMathAnal,FengGuanLowengrub2018_JSciComput,GompperSchick1990_PRL,KraitzmanPromislow2018_SIAMJMathAnal}.
These refined continuum formulations are in principle capable of capturing
topological changes of a lipid bilayer membrane, such as membrane fusion and fission. 
However, no simulations of membrane fusion or fission based on these refined formulations are available in the literature (to our knowledge), possibly
due to the numerical challenges to efficiently and accurately resolve structures on the scale of membrane thickness.

Changes in topology of bilayer membranes, as occur in bilayer membrane fusion, pore formation and protein insertion, for example,
involve the introduction of a hydrophobic fissure
in the normally intact monolayer surface. 
Due to the relatively large tension of a hydrocarbon-water interface, the energy of a hydrophobic fissure 
can dominate the membrane's elastic energy,
making it necessary to also take into account local interactions at the molecular level~\cite{FrRoPi17,Discher967}. 
Moreover, in many subcellular structures, membrane energies are dominant
in high curvature regions only a few lipids wide~\cite{Huetal08,Yoo2013}. 

Based on these observations, we focus on topological changes with mesoscopic interactions 
in a semi-continuum framework where the lipids are coarse-grained into amphiphilic Janus-type particles while their
interactions with each other and the solvents are described at a continuum level.
This hybrid approach provides a bridge from microscopic molecular formulation to macroscopic continuum description of a lipid bilayer membrane.
Furthermore, the continuum limit of our hybrid mesoscopic model may facilitate efficient numerical algorithms for simulating fusion/fission of lipid bilayer membranes
of physically relevant membrane size and dynamic duration.

Modern molecular dynamics (MD) simulators (such as MARTINI~\cite{Marrink07} and\break LAMMPS~\cite{plimpton2007lammps}) 
have the advantage in that they resolve all relevant  molecular details,
and have been widely employed to simulate
fully atomistic or coarse-grained lipid bilayer membrane based on pairwise interactions~\cite{Cooke05,Fu16,JuVa17,KaVa18,Marrink07,Atzberger}.
Traditionally, MD numerical methods use explicit fluid particles such as coarse-grained water molecules 
and pairwise Lennard-Jones interactions.
There is the disadvantage, though, that an enormous number of water molecules and long computation time are needed in MD simulations, and
it remains a great challenge to compute the hydrodynamic interactions of the lipid membrane at micron size for durations
long enough to make physical predictions. 
One way to mitigate long computation time is to compute hydrodynamic interactions using an implicit solvent and
Stokesian dynamics \cite{BrBo88}.


In the present work, we propose a novel approach to lipid-lipid interactions
called the hydrophobic attraction domain functional (HADF).
Let $\Omega$ be an open, exterior domain in $\mathbb{R}^n$ 
representing water surrounding a collection of amphiphilic particles, e.g. lipids.
For dispersed particles, the energy associated with hydrophobic interfaces
behaves as a surface energy. When nearby, particles decrease their energy by
aggregating and sequestering their hydrophobic interfaces from water. 
These interactions are well-described by the 
Ginzburg-Landau-type domain functional
\begin{equation}
\label{eq:main}
\Phi(\Omega,f) = \gamma  \min_{u \in \mathcal{A}} I[u],
 \end{equation}
where
\begin{equation}
\label{eq:main2}
 I[u]  = \int_{\Omega} \rho |\nabla u|^2 + \rho^{-1} u^2 \,dx.
\end{equation}
Here $\mathcal{A} = \{u \in W^{1,2}(\Omega) : u = f \text{ on } \Sigma\}$ 
is the admissible class and $f$ with range $[0,1]$ is the hydrophobicity label for the 
water-particle interface $\Sigma = \partial \Omega.$
The parameter $\gamma > 0$ is interfacial tension. 
Its value in bilayers has been widely investigated in both numerical and 
theoretical studies \cite{EVANS03,GarciaSaez,Nagle17,Petelska2012}.
For a Lipschitz domain $\Omega$ and  for $f$ the trace of a function in $W^{1,2}(\Omega)$ \cite{EVANS},
the existence of a unique minimizer to \cref{eq:main} is a straightforward consequence of the closest point 
theorem~\cite{LAX}.

The scalar function $u$ of \cref{eq:main} models disruption in the hydrogen 
bonding network \cite{ErLjCl89,MeMeFe17}. For a point $x \in \Sigma$ 
representing a hydrophobic interface, water mobility is restricted and there $u(x) = f(x) = 1.$
Conversely, $u(x) = f(x) = 0$ at a point $x \in \Sigma$ representing a hydrophilic interface where water mobility
is unrestricted. 
In the water region, $u$ in \cref{eq:main} is a solution to the boundary value problem (BVP) of the screened Laplace equation:
\begin{equation}
        \label{SL}
\begin{cases}
 -\rho^2 \Delta u + u =0 & \text{ in } \Omega,\\
  u(x) = f(x)  &\text{ on }\Sigma, \\
   u(x) \to 0, &\text{ as } x \to \infty.
\end{cases}
\end{equation}
Solutions $u$ of \cref{SL} have a boundary layer of thickness $\rho>0.$
Thus disruption in hydrogen bonding modeled by \cref{SL} extends into the bulk
with characteristic distance $\rho$ \cite{DuMa16,MeRoIs06}.

\begin{figure}[t]
\begin{center}
\includegraphics[width=0.32\textwidth]{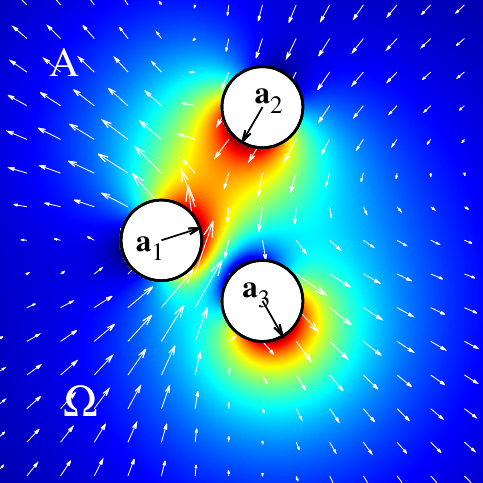}
\includegraphics[width=0.32\textwidth]{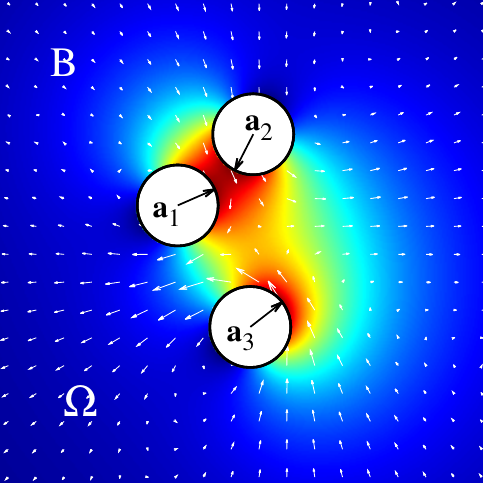}
\includegraphics[width=0.32\textwidth]{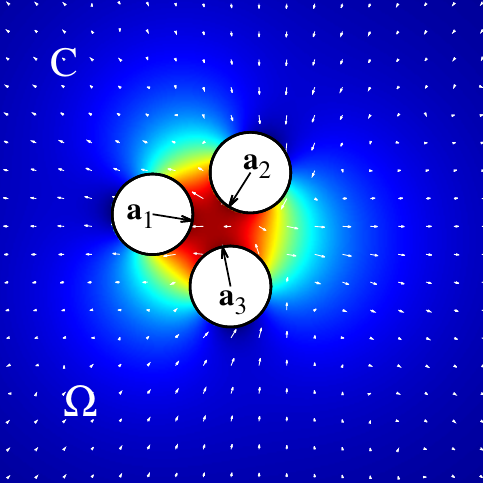}
\end{center}
  \caption{
  The figure illustrates hydrophobic attraction of amphiphiles in a zero-Reynolds number fluid.   
  The black arrows are the directors $\mathbf{d}_i$ of the particles centered at $\mathbf{a}_i$, $i = 1,2,3$.
  The white arrows are the fluid velocity. 
  The color map has dark blue is for $u = 0$ (lipid heads) and dark red is for the hydrophobic interface
  $u = 1$ (lipid tails). Activity, shown as pseudo-color shading in the figure, extends from the 
  hydrophobic interface into the bulk. 
  Going from panels A to B, particles $\mathbf{a}_1$ and $\mathbf{a}_2$ come together 
  and $\mathbf{a}_3$ rotates counterclockwise. 
  In panel C, the particles sequester activity 
  to a single hydrophobic core. 
  }
    \label{figure1}
\end{figure}

The hydrophobic force is the first variation of the functional $\Phi$ 
with respect to the shape of the domain $\Omega.$ 
The challenge in the present work is to  compute the hydrophobic force
between several bodies of arbitrary shape and configuration. 
\Cref{sec:variation} carries out the variation for rigid body motions, and this reduces to a set of 
boundary integrals for the hydrophobic forces and torques. 
For simulations, we utilize a boundary condition $f$ representing surface portions of lipid tail and lipid head, 
and we adopt an excluded volume repulsion to avoid particle collisions in the many-body simulations.
As an illustration, \cref{figure1} shows
the self-assembly process for three Janus-type particles in a viscous fluid.

 An important feature of the model is that the potential $\Phi$ and its intermolecular forces and torques, in contrast to that of coarse-grained theories, do not arise from any pairwise potentials (see \cref{sec:non-pairwise}).
To leading order, the attraction between particle pairs predicted by \cref{eq:main} 
is in accord with experimental force-distance curves \cite{ErLjCl89, Lietal05, MeRoIs06}.
The functional \cref{eq:main}, however, requires modification for account for sub-nanometer
force oscillations observed in experiment,
e.g. through the inclusion of higher order terms.
Nevertheless, the HADF captures the essential features of amphiphile self-assembly,
and the variational calculations and numerical methods generalize to more
complicated domain functionals.

An essential principle for molecular or particle based approaches is to 
ensure that the total free energy accounting 
for lipid-lipid and lipid-water interactions 
gives rise to an equivalent elastic characterization of membranes as 
determined by experimental measurements \cite{TeDe18,VeBrPa15}.
\Cref{sec:elasticity} examines the elasticity of bilayer particle configurations.
We obtain physical quantities such as bending modulus, tilt modulus and stretching modulus 
by setting up corresponding  equilibrium simulations from continuum theory \cite{KoNa15,Naetal15,TerziDeserno17}.

Section \ref{sec:dynamics} formulates the mobility problem to calculate hydrodynamics from the hydrophobic stress.
The dynamics for many-particle simulations yield physically reasonable 
time scales and configurations.
For example, we can track the particle dynamics over the nanosecond range needed for rapid particle self assembly, up to the microseconds range where bilayer and micelle shapes evolve over a slower time scale
\cite{Reetal12, Rietal18}.

Calculating the particle dynamics requires
rapid, on-the-fly solution of \cref{SL}. In \cref{sec:bie}, we 
present a new SKIE formulation for the boundary value problem \cref{SL}, derived
from a representation of the solution in which
the unknowns are only on the boundary $\Sigma$. In \cref{sec:numerical_methods}, we describe an approach to applying a recently developed QBX-FMM
scheme for discretizing the SKIE accurately and adaptively, solving
the resulting linear system and evaluating the desired physical quantities afterwards accurately and rapidly. 
The resulting scheme has linear
complexity with an optimal number of unknowns for the simulation of particle dynamics at each time step.
To compare the computational cost against MD simulations, even solvent free coarse-grained models have at least $O(N^2)$ complexity in the  number of particle $N$ \cite{Cooke05,noguchi11}.


\section{Intermolecular Forces and Torques}
\label{sec:variation}
\begin{figure}
\begin{center}
\includegraphics[width=0.66\textwidth]{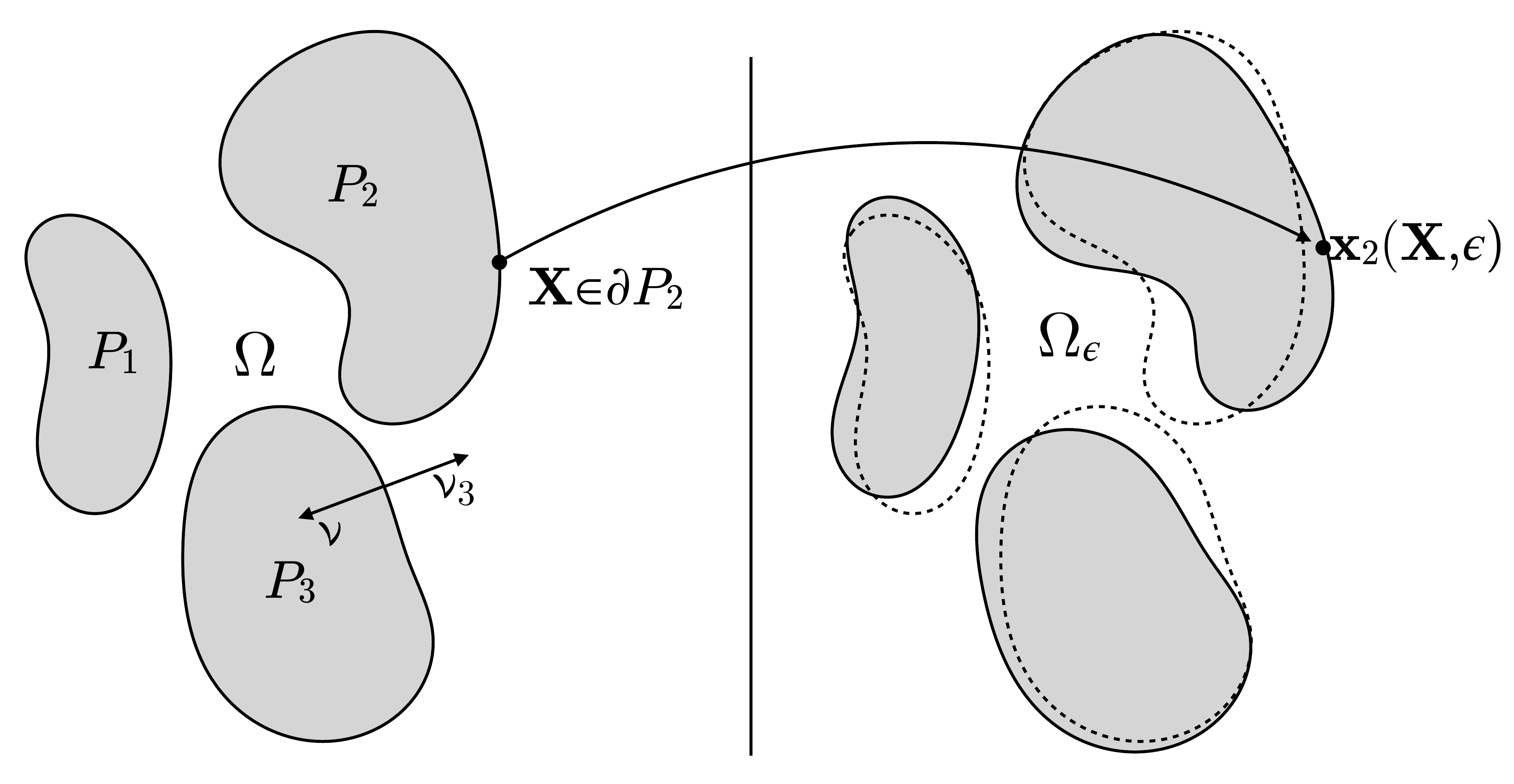}
\end{center}
\caption{\label{figure2} On the left are three particles $P_1,$ $P_2$ and $P_3$
forming the exterior domain $\Omega$. Rigid transformation of each particle (right panel) leads to the 
perturbed exterior domain $\Omega_{\epsilon}$ and changes in the relative positions of the material
label $f_{\epsilon}$, resulting in variations in the hydrophobic potential $\Phi$. }
\end{figure}

We calculate the first variation of $\Phi$ with respect to rigid body deformations \cite{BaWa15,Schiffer1954}.
Consider $N$-many, rigid particles represented by 
disjoint, bounded, closed regions $P_1, P_2, \dots, P_N$  in  $\mathbb{R}^{n},$  $n = 2,3$. 
The water region (the exterior domain) and particle-water interface are 
\begin{equation}
\Omega = \mathbb{R}^{n} \setminus \bigcup_{i=1}^N P_i,\quad  \Sigma = \bigcup_{i=1}^N \partial P_i,
\end{equation}
respectively. Throughout, $\nu$ denotes the unit outward normal to $\Omega$, and $\nu_i$ 
denotes the unit outward normal to $P_i$. 
Note that $\nu$ and $\nu_i$ have opposite orientation, as illustrated in  Figure~\ref{figure2}.
Suppose that $u$ is the solution to the screened Laplace BVP \cref{SL} with the material label $f$. 
Then the force $\mathbf{F}_i$ and torque $\tau_i^0$ acting on 
particle $P_i$ are
\begin{equation}
\label{eq:force_torque}
\mathbf{F}_i
= \int_{\partial P_i}\mathbf{T}\cdot \nu_i\, dS,\quad
\tau_i^0 
= \int_{\partial P_i}\mathbf{r}_0 \times (\mathbf{T}\cdot \nu_i)\, dS,
\end{equation}
where  
\begin{equation}
\label{eq:stress}
\mathbf{T}
= \gamma\rho^{-1}u^2 \mathbf{I} + 2\rho\gamma (\tfrac{1}{2}|\nabla u|^2 \mathbf{I} - \nabla u\otimes \nabla u),
\end{equation}
is the hydrophobic stress and $\mathbf{r}_0$ is the position vector relative to the origin $\mathbf{0}$.
To ensure that \cref{eq:force_torque} is well-defined 
and to guarantee differentiability of the domain functional,  
we that $\Omega$ is a $C^{2,\alpha}$ domain and that 
$f = \tilde{f}$ on $\Sigma$ for some $\tilde{f} \in C^{2,\alpha}(\overline{\Omega})$. 

To compare $\Phi(\Omega,f)$ against that of competing domains, 
consider a one-parameter family of rigid transformations 
\begin{equation}
\label{eq:rigid}
\mathbf{x}_i(\mathbf{X},\epsilon) = \mathbf{c}_i(\epsilon) + \mathbf{Q}_i(\epsilon)\mathbf{X}, 
\end{equation}
parametrized by $\epsilon \in \mathbb{R}$. 
The vector $\mathbf{c}_i(\epsilon)$ and tensor $\mathbf{Q}_i(\epsilon)$ give the displacement and rotation 
of the particle $P_i$, $i = 1,\dots, N$, 
relative to the origin. They satisfy $\mathbf{c}_i(0) = \mathbf{0}$ and $\mathbf{Q}_i(0) = I$ so that $\mathbf{x}_i(\mathbf{X},0)$ 
is the identity transformation; $\mathbf{Q}_i(\epsilon)\mathbf{Q}_i^T(\epsilon) = I$ for all $\epsilon$. 
The distance between $P_i$ and $P_j$ is positive whenever $i \neq j$. Therefore, for $\epsilon$ in an open interval about $0$, let 
\begin{equation}
\label{transportDomain}
\Omega_\epsilon = \mathbb{R}^n \setminus \bigcup_{i=1}^N \mathbf{x}_i(P_i,\epsilon),\quad \Sigma_\epsilon = \partial \Omega_{\epsilon},\quad
f_\epsilon(\mathbf{x}_i(\mathbf{X},\epsilon)) = f(\mathbf{X}), \quad \mathbf{X} \in \partial P_i.
\end{equation}
Finally, let $u_\epsilon(x)$ be the one-parameter family of solutions to the perturbed boundary value problem of screened Laplace equation 
\begin{equation}
\label{eq:transportSL}
-\rho^2 \Delta u_\epsilon + u_\epsilon = 0 \text{ in } \Omega_\epsilon,\quad  
u_\epsilon = f_\epsilon \text{ in } \Sigma_\epsilon,\quad
u_{\epsilon} \to 0 \text{ as } x \to \infty.
\end{equation}
The domain $\Omega_\epsilon$ and boundary $\Sigma_\epsilon$ are the water region and water-molecule interface after 
transforming
each particle according to its rigid motion \cref{eq:rigid} (see \cref{figure2}).
For $x \in \Omega$, let 
\begin{equation*}
\dot{u}(x) = \frac{d}{d\epsilon}u_\epsilon(x)\Big|_{\epsilon=0},
\end{equation*}
and extend $\dot{u}$ continuously to $\overline{\Omega}$. 
Due to \cref{transportDomain}, we have the transport identity 
\begin{equation}
\label{eq:transport}
\dot{u} + \nabla u  \cdot \dot{\mathbf{x}} = 0 \mbox{ on } \Sigma,
\end{equation}
where $\dot{\mathbf{x}}(\mathbf{X}) = \frac{d\mathbf{x}_i}{d\epsilon}(\mathbf{X},0)$ whenever $\mathbf{X} \in P_i$. 
(Note, however, that the values of $\dot{u}$ in $\Omega$ are determined by the BVP \cref{eq:transportSL}, 
and therefore do not generally satisfy this transport relation.) 

Applying  the Reynolds transport theorem~\cite{LEAL}, we obtain
\begin{equation}
\begin{aligned}
&\frac{d}{d\epsilon}\Phi(\Omega_\epsilon,f_\epsilon)\Big|_{\epsilon=0} 
=  \gamma \frac{d}{d\epsilon} \left(\int_{\Omega_\epsilon} \rho|\nabla u_\epsilon|^2 + \rho^{-1} u_\epsilon^2 \,dx\right)\Big|_{\epsilon=0} \\
&=  \gamma \int_{\Omega} 2\rho \nabla u \cdot \nabla \dot{u} + 2\rho^{-1} u \dot{u} \,dx 
+ \gamma \int_{\Sigma} \left(\rho|\nabla u|^2 + \rho^{-1} u^2\right) \dot{\mathbf{x}}\cdot \nu \,dS.
\end{aligned}
\end{equation}
Integration by parts then gives
\begin{align}
\label{eq:int}
\frac{d}{d\epsilon}\Phi(\Omega_{\epsilon},f_{\epsilon})\Big|_{{\epsilon}=0}= 
\gamma \int_{\Sigma} \left(\rho|\nabla u|^2 + \rho^{-1} u^2\right) \dot{\mathbf{x}}\cdot \nu - 2\rho\nabla u \cdot \nu \dot{u}\,dS.
\end{align}
Due to the minimality condition $-\rho^2\Delta u + u = 0$, the interior values of $\dot{u}$  do not enter \cref{eq:int}. 
Based on \cref{eq:transport} and the fact that $\nu$ and $\nu_i$ have opposite orientation, 
\begin{align}
\notag
\frac{d}{d\epsilon}\Phi(\Omega_{\epsilon},f_{\epsilon})\Big|_{{\epsilon}=0}&=
\gamma\sum_{i=1}^N \int_{\partial P_i} -\left(\rho|\nabla u|^2 + \rho^{-1}u^2\right) \nu_i \cdot \dot{\mathbf{x}}_i(0)
 + 2\rho \nabla u \cdot \nu_i  \nabla u \cdot \dot{\mathbf{x}}_i(0) \,dS\\
\notag
 &= \gamma\sum_{i=1}^N \int_{\partial P_i} \dot{\mathbf{x}}_i(0) \cdot \left[-\rho^{-1}u^2 \mathbf{I} + 2\rho(\nabla u \otimes \nabla u
 - \tfrac{1}{2}|\nabla u|^2 \mathbf{I})\right]\cdot \nu_i \,dS\\
\notag
 &= -\sum_{i=1}^N \int_{\partial P_i} (\dot{\mathbf{c}}_i(0)  + \dot{\mathbf{Q}}_i(0) \mathbf{r}_0) \cdot \mathbf{T} \cdot \nu_i\,dS\\
\label{var2force}
&= -\sum_{i=1}^N \left(\dot{\mathbf{c}}_i(0) \cdot \mathbf{F}_i + \mathbf{w}_i \cdot \tau_i^0\right),
\end{align}
where $\mathbf{w}_i = \langle w_{1}^i,w_{2}^i,w_{3}^i\rangle$ 
is the axial vector for the skew symmetric tensor $\dot{\mathbf{Q}}_i(0)$. 
In the second to last equation, the minus sign makes the force act in the negative direction
of the potential gradient. 
This establishes \cref{eq:force_torque} and \cref{eq:stress}.

In the formulation \cref{eq:rigid}, the rigid motions are independent.
Consider the case when the rigid motions are uniform, that is, 
$\mathbf{c}_i(\epsilon) = \mathbf{c}(\epsilon)$and $\mathbf{Q}_i(\epsilon) = \mathbf{Q}(\epsilon)$ for all $i = 1,\dots, N$. Then 
the solution to the perturbed BVP (\ref{transportDomain}, \ref{eq:transportSL}) satisfies
\begin{equation}
\label{puretransport}
u_\epsilon(\mathbf{c}(\epsilon) + \mathbf{Q}(\epsilon)\mathbf{X}) = u(\mathbf{X}).
\end{equation}
It follows that $\Phi(\Omega_\epsilon,f_\epsilon) = \Phi(\Omega,f)$ for all $\epsilon$ and, by \cref{var2force}, that
\begin{equation}
\label{eq:constant}
\sum_{i=1}^N \dot{\mathbf{c}}(0) \cdot \mathbf{F}_i + \mathbf{w} \cdot \tau_i^0
= -\frac{d}{d\epsilon}\Phi(\Omega_\epsilon,f_\epsilon)\Big|_{\epsilon=0} 
= -\frac{d}{d\epsilon}\Phi(\Omega,f)\Big|_{\epsilon=0} = 0.
\end{equation}
Here,  $\mathbf{w} = \langle w_{1},w_{2},w_{3}\rangle$ is the axial vector for $\dot{\mathbf{Q}}(0)$. 
Since $\dot{\mathbf{c}}(0)$ and $\mathbf{w}$ are arbitrary, we have 
\begin{equation}
\label{eq:net}
\sum_{i=1}^N \mathbf{F}_i = \mathbf{0},\quad \sum_{i=1}^N \tau_i^0 = \mathbf{0}.
\end{equation}
In other words, the net hydrophobic interaction 
is force and torque free. 

 
 \subsection{Simulations}
 \label{Simulations}
 For the simulations in this paper, 
the $P_1, \dots, P_N$ are two-dimensional Janus-type particles.
The direction vector $\mathbf{d}_i = \langle \cos \theta_i, \sin \theta_i \rangle$  
specifies  orientation and $\mathbf{a}_i$ is the 
particle position (e.g. the center of mass, \cref{figure1}). The particle shapes are ellipses with  
semi-major and semi-minor axes $a_i$ and $b_i,$ respectively. 
In the case of lipids, $2a_i$ represent lipid length and 
major axis is parallel to the director and hydrocarbon tail.
 
The material label for the Janus-type particle takes the form
\begin{equation}
\label{JanusBC}
f(\mathbf{x}) = 1 - \sin^p(\theta), \quad \mathbf{x} \in \partial P_i,
\end{equation}
where $\theta$ is the angle formed by $\mathbf{x}-\mathbf{a}_i$ and $\mathbf{d}_i.$
Accordingly, there is a smooth transition in hydrophobicity across the particle \cite{MacCallumTieleman}, 
with the boundary portion in the direction $\mathbf{d}_i$ modeling a hydrophobic
tail and the opposite boundary portion modeling a hydrophilic head. 
The size of the hydrophobic region grows with the even integer parameter $p.$ 
Finally, 
\begin{equation}
\label{eq:scalar-torque}
\tau_i = \tau_i^0 - \mathbf{a}_i \times \mathbf{F}_i
\end{equation}
is the two-dimensional (scalar) torque about the position $\mathbf{a}_i$. 

For small but fixed separations between particles, 
our numerical scheme accurately
resolve the field $u$  without an undue cost increase due to refinement;
we postpone a detailed discussion of the method and involved cost to Section~\ref{sec:numerical_methods}.
Dynamically, the forces \cref{eq:force_torque} bring the coarse-grained lipid particles into contact.
An excluded volume repulsion prevents near-contact between particles \cite{ParsegianNinham}. 
For two circular particles, the interaction is 
\begin{equation}
\label{repul}
\mathbf{F}^{\mathrm{rep}}_{ij} = c_0 \frac{q}{(|\mathbf{a}_i - \mathbf{a}_j|-(b_i+b_j))^{q+1}}\frac{\mathbf{a}_i - \mathbf{a}_j}{|\mathbf{a}_i - \mathbf{a}_j|}, \quad  i\neq j.
\end{equation}
We fix the order $q = 3$ ($q = 4$ in three-dimensions)
and use the parameter $c_0$ to control the strength of repulsion.
For ellipses of eccentricity close to zero, we approximate the excluded volume
repulsion using three circular particles placed along the major axes,
as described in Supplementary Material, Section S1. In the sequel, 
\begin{equation}
\label{eq:ellipse_rep}
\mathbf{F}^{\mathrm{rep}}_i,\quad \tau_i^{\mathrm{rep}}, \quad
\Phi_{\mathrm{rep}},
\end{equation} 
denote the excluded volume force, torque and repulsion potential, respectively.
The total potential that includes hydrophobic attraction and 
steric repulsion is 
\begin{equation}
\label{total_en}
\Phi_{\mathrm{Total}} = \Phi + \Phi_{\mathrm{rep}}.
\end{equation}

For the simulations, 
we assume translation invariance in the $z$-direction. Figures \ref{fig:stretching} and \ref{figure4} give values 
in $\KBT$ per length 
since the two-dimensional
simulations are for the cross-section
of a three-dimensional bilayer. All 
other physical parameters correspond to their
usual three-dimensional value.

We use $2a_i = 2.5$ nm as a representative phospholipid length \cite{Boal},
the screening length $\rho = 2.5$ nm \cite{ErLjCl89,Lietal05,Parsegian,Israelachvili80,TerziDeserno17}, and
$c_0=0.5$ pN nm$^4$ for the inter-particle repulsion.
Bilayers containing different single pure components give 
various interfacial tension $\gamma$
values which are within the range of $0.7$ -- $5.3$ pN nm$^{-1}$  \cite{KUZMIN2005, Petelska2012}. We find that the mechanical moduli calculated 
from our simulation data are in good agreement with 
results in the experimental literature when the interfacial tension
$\gamma=4.1$ pN nm$^{-1}.$ Coincidentally, this value corresponds to a 
specific lipid composition DPoPC:SM:Chol in bilayer membrane \cite[Table 1]{GarciaSaez}.

%
Our experiences show that the computational cost to maintain the same
numerical accuracy in solving the boundary value problem \eqref{SL} grows only moderately when going from circular to elliptical model particles. For instance,
ellipses with $a_i/b_i = 3$ require 60 \% more grid points than for $a_i/b_i=1.$
At the same time, ellipses afford flexibility in terms of dimensions that determine physical properties of bilayer. 
However, we remark that rather than representing
a physical lipid or collection of lipids, the model
particle discretizes the mean lipid position and orientation
but without the mesh associated with finite element methods, for example 
\cite{Bartels,Ryham16}.
Similarly, the gap region between 
neighboring particles indicates a hydrophobic zone and not an intervening water.

\section{Bilayer Elasticity}
\label{sec:elasticity} 
\begin{figure}[t]
\begin{center}
\includegraphics[width=0.99\textwidth]{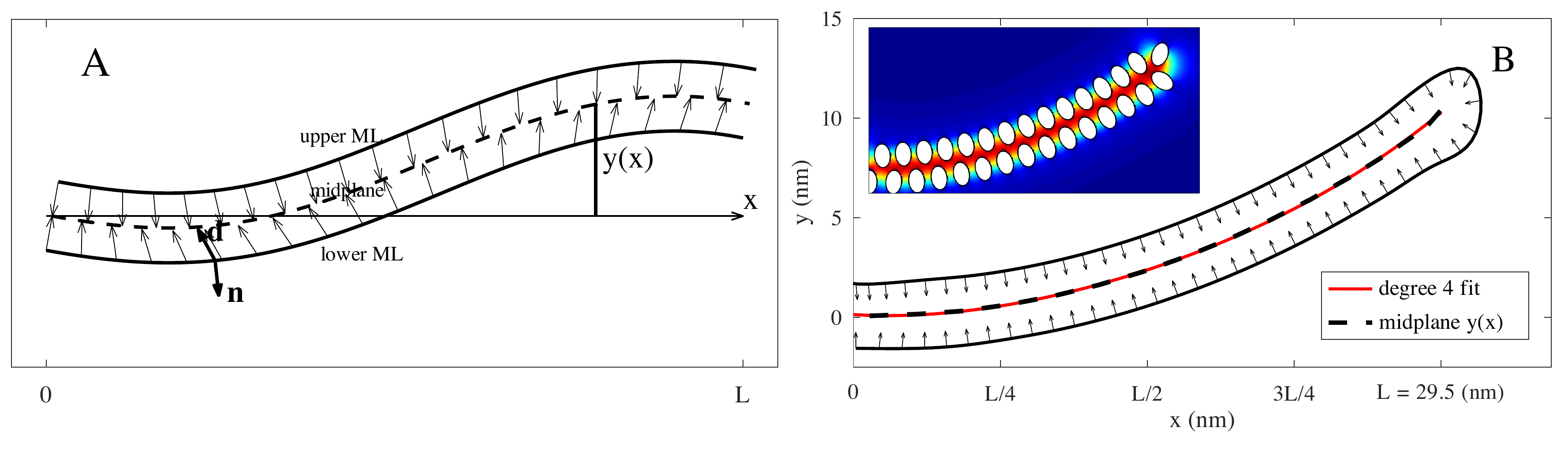}
\end{center}
\caption{Panel A depicts a bilayer, with its upper and lower monolayers
(solid curves) and midplane (dashed curve). 
The director field $\mathbf{d}$ points in the direction of the 
lipid tails and the normal $\mathbf{n}$ points outward from
bilayer core. When deformations and tilt are small 
($|y'| \ll 1$, $|\mathbf{d} \times \mathbf{n}| \ll 1$)
energy is quantified in terms of the midplane has the height function $y(x).$ 
Panel B has the equilibrium configuration
Janus-type amphiphilic particles with uniform loading and a clamped boundary 
condition.
Fitting \cref{eq:zero_tilt_sol} to the midplane (dashed curve) yields 
the bending modulus $\KB = 8.51$ $\KBT.$
In Panel B (inset), enumerating counterclockwise from the bottom left, 
particles satisfy 
$x_1 = 0,$ $x_{16} = 1.5$, $\theta_1 = \pi/2,$ $\theta_{16} = -\pi/2$ 
and  $y_1=0.$ The parameters are 
$\gamma = 4.1$ pN nm$^{-1}$, 
$\rho   = 2.5$ nm, 
$c_0    = 0.5$ pN nm$^{-4}$ and
$k = 0.0116$ pN nm$^{-2}.$
The ellipses have $a_i = 1.25$ nm and $b_i = 0.8$ nm.
 }
\label{fig:bend}
\end{figure}

We compare our two-dimensional 
equilibrium configurations to those found in membrane continuum mechanics. 
In large particle number HADF simulations, particles bring opposing hydrophobic 
regions into contact, forming two abutting monolayers of a bilayer.  
Continuum theory describes monolayers using 
a director field $\mathbf{d}$ to track lipid orientations, 
along with a field $\mathbf{n}$ given by the monolayer
surface normal (Figure \ref{fig:bend}A),    
and  quantifies monolayer energy using a 
Helfrich Hamiltonian 
\begin{equation}
\label{eq:Helfrich}
\int_{\mathcal{C}} 
\tfrac{1}{2}\KB\left[ \left( \Div \mathbf{d} + k_0\right)^2 - k_0^2\right] 
+ \tfrac{1}{2}\KTH |\mathbf{d} \times \mathbf{n}|^2 
\,ds. 
\end{equation}
The curve $\mathcal{C}$ tracks the cross-section of the monolayer neutral
surface. The integrand in \cref{eq:Helfrich} contains
the splay distortion $\Div \mathbf{d}$ 
with bending modulus $\KB$, and the tilt deformation $\mathbf{d} \times \mathbf{n}$ 
with tilt modulus $\KTH$ \cite{Nagle17-2}.  
The parameter $k_0$ is spontaneous curvature \cite{RoLi15,Kozlov2007,FriedSeguin15}.
Since we are assuming translational invariance in the $z$-direction, the twist $\Curl \mathbf{d}$
and saddle-splay $\det \mathsf{D}$ distortions are absent from \cite{Helfrich73,TerziDeserno17},
and \cref{eq:Helfrich} behaves as an energy density per length.

Consider a planar bilayer subject to a uniform vertical load. The bilayer is clamped and horizontal 
at one end and the restoring force of bending in the free part of the bilayer opposes the load.
Taking $\mathbf{d}$ parallel to $\mathbf{n}$ and assuming a small deformation
gives the appropriate functional
\begin{equation}
\label{eq:bend_functional}
\int_0^L \KB (y'')^2 - ky\,dx, \quad y(0) = y'(0),
\end{equation}
where $y(x)$ is the height function for the bilayer midplane (Figure \ref{fig:bend}A, dashed curve),
$\Div \mathbf{d} = \pm y''$ and $k$ is the load strength.
The summation of the monolayer energies \cref{eq:Helfrich} with opposite normals leads to the cancelation of the spontaneous curvature terms in \cref{eq:bend_functional}.

Minimizers of \cref{eq:bend_functional} satisfy
the boundary value problem 
\begin{equation}
2\KB y^{(4)} = k,\quad y(0) = y'(0) = y''(L) = y^{(3)}(L) = 0.
\end{equation}
We find the solution 
\begin{equation}
\label{eq:zero_tilt_sol}
y(x) = \frac{kL^4}{2\KB}\left[ \frac{1}{24}\left(\frac{x}{L}\right)^4 - \frac{1}{12}\left(\frac{x}{L}\right)^3 + \frac{1}{4}\left(\frac{x}{L}\right)^2 \right].
\end{equation}
Thus, we can determine $\KB$ from   curves of the form \cref{eq:zero_tilt_sol}
whenever $L$ and $k$ are given.

The inset in Figure \ref{fig:bend}B shows a HADF equilibrium 
configuration used to determine $\KB.$
The $N = 30$ particles minimize the modified functional 
\begin{equation}
\label{eq:HAFload}
\Phi_{\mathrm{Total}} - \sum_{i=1}^N \tilde k y_i, 
\end{equation}
where the $\tilde k = Lk/N$ is the discrete load strength coming from quadrature of the integral 
\cref{eq:bend_functional} with $N$ many particles. To achieve minimality,
the particles start in the shape of a flat bilayer, and then migrate upward 
following steepest gradient for \cref{eq:HAFload}.

The main figure in Figure \ref{fig:bend}B depicts the monolayer neutral  surface
(solid curve), midplane (dashed curve) and the lipid directors interpolated  
from the discrete particle positions and orientations (of the inset). 
The directors are everywhere normal to the neutral surface and the 
deformations  are small. 
This justifies applying the zero-tilt, small-deformation solution \cref{eq:zero_tilt_sol}.  
Fitting a 4th degree polynomial to the midplane curve (Figure \ref{fig:bend}B, red curve)
supplies the coefficient $(kL^4/2\KB)$ of \cref{eq:zero_tilt_sol}. 
Combining the coefficient with simulation values for $L$ and $k$ (Figure \ref{fig:bend}B, caption)
yields $\KB = 8.51$ $\KBT.$ 
This value for the bending modulus is for ellipses using $p = 6$ in the hydrophobicity boundary condition \cref{JanusBC}.
To assess how bilayer rigidity depends on the material label, we considered 
the energy minimization with $p = 2,$ which gave  $\KB =  13.54$ $\KBT.$
We conclude that under HADF, particle configurations behave like an elastic material.
The associated bending modulus grows with symmetry in the hydrophobic surface label, e.g.
$\KB$ was largest for $p = 2$ where the label is symmetric across $\theta = \pi/2.$

\begin{figure}
\begin{center}
\includegraphics[width=0.99\textwidth]{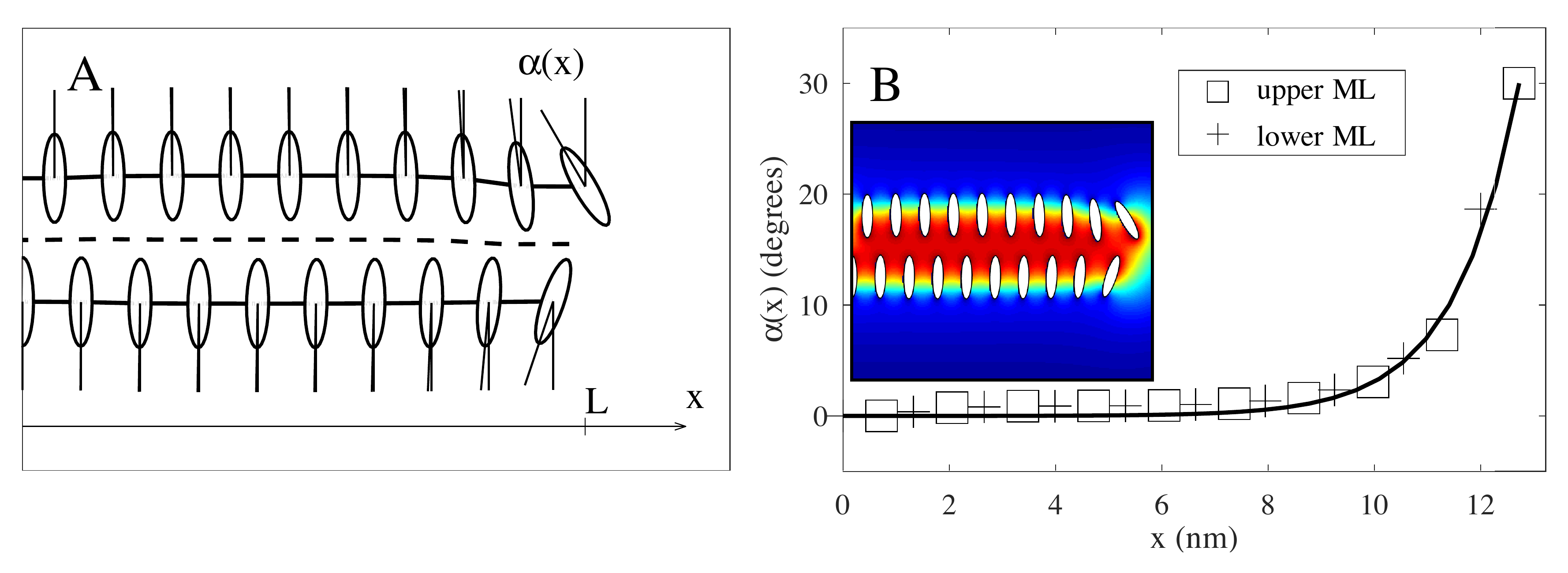}
\end{center}
\caption{The monolayers in panel A are flat and tilt is nonzero. 
  The horizontal coordinate $x$ runs between 0 and the length of bilayer, $L,$  
  and the function $\alpha(x)$ measure the angle between $-\mathbf{d}$ and the unit normal.
  The boundary conditions for \cref{eq:tilt_energy} are $\alpha_0 = 0^{\circ}$ and $\alpha_1 = 30^{\circ}.$ 
  In panel B, the $+$ and $\square$ symbols are the angles at the particle centers.
  The parameters are 
$\gamma = 4.1$ pN nm$^{-1}$,
$\rho   = 2.5$ nm,
  $c_0    = 0.5$ pN nm$^{-4}$, and
  $k = 0.0116$ pN nm$^{-2}.$
The ellipses have $a_i = 1.25$ nm and $b_i = 0.3125$ nm.
  }
    \label{fig:tilt}
\end{figure}

Now we consider a flat monolayer 
with nonzero tilt (Figure \ref{fig:tilt}A). The splay distortion comes from changes in the angle $\alpha(x)$ 
between the director $\mathbf{d}$ and the vertical. 
For small angles, the monolayer energy \cref{eq:Helfrich} becomes 
\begin{equation}
\label{eq:tilt_energy}
\int_0^L \tfrac{1}{2}\KB (\alpha')^2 + \tfrac{1}{2}\KTH \alpha^2 \,dx,\quad
\alpha(0) = \alpha_0,\quad \alpha(L) = \alpha_1.
\end{equation}
Note that we have left off the null-Lagrangian term $\KB k_0 \alpha'$
from this expression since $\KB,$  $k_0$ 
and the boundary data $\alpha_0$ and $\alpha_1$ are constants.
Assuming $\alpha(0) = 0,$ minimizers of \cref{eq:tilt_energy} take the form 
\begin{equation}
\label{eq:tiltfit}
\alpha(x) = \alpha_1\frac{\sinh(x/\kappa)}{\sinh(L/\kappa)},
\end{equation}
where $\lambda = \sqrt{\KB/\KTH}$ is the tilt decay length \cite{KUZMIN2005}. 
Figure \ref{fig:tilt}B shows the data (plusses and squares) for 
the HADF equilibrium configuration with fixed endpoint angles.
The solid curve fits \cref{eq:tiltfit} to the 
angle data for the value $\lambda = 1.2$ nm. 
This value is consistent with experimental and theoretical measurements 
of the bending and tilt moduli \cite{Nagle17-2, KoNa15}.

In HADF, tilt dissipation is a consequence of repulsion between rod-like particles.
The ellipses in Figure \ref{fig:tilt} are elongated and have $a_i/b_i = 4.$ 
When the particles are more circular ($a_i/b_i \sim 1$), the bulk particles 
ignore endpoint orientations and the angle function $\alpha(x)$ is non-monotonic in $x.$

\begin{figure}
\begin{center}
\includegraphics[width=\textwidth]{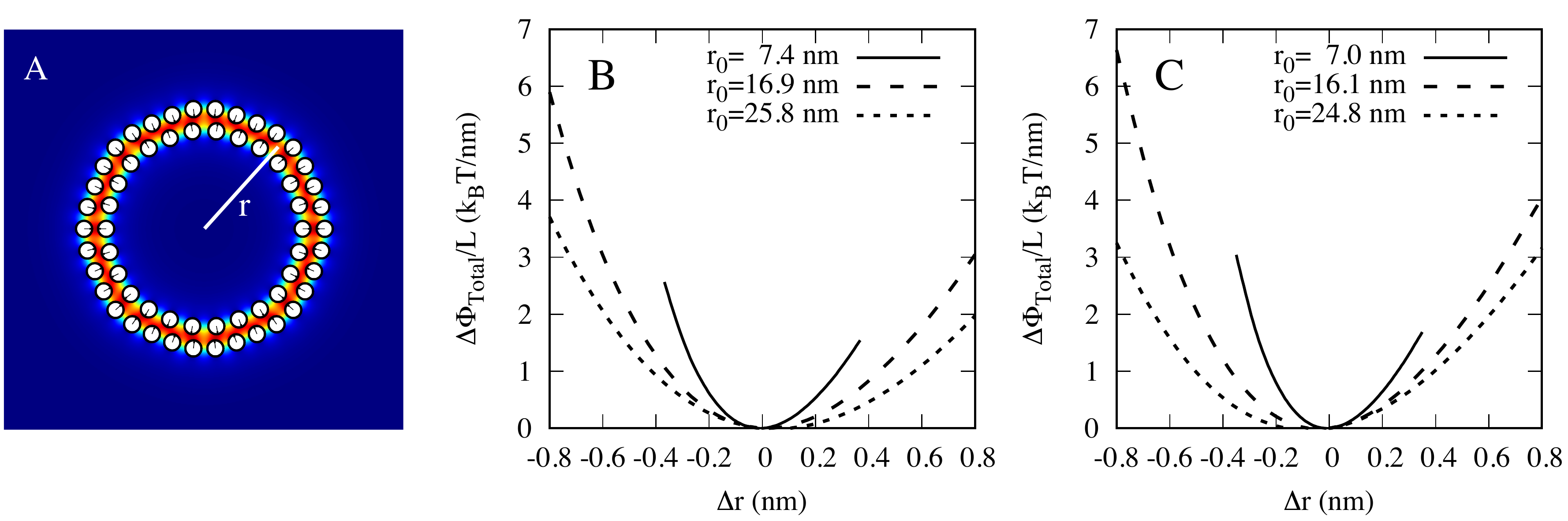}
\end{center}
\caption{
  Panel A shows the activity field for the cross-section of a cylindrical bilayer.
  The radius measures from the center to the midplane.
  Panels B and C plot the
  change in energy $\Delta \Phi_{\mathrm{Total}} = \Phi_{\mathrm{Total}}(r) - \Phi_{\mathrm{Total}}(r_0)$ per length 
  of the cylinder under stretching and compression: $\Delta r = r - r_0.$
  Equilibrium radius $r_0$ decreases with $c_0;$
  the curves are for $c_0 = 0.5$ pN nm$^4$
  in panel B,    and $c_0 = 0.25$ pN nm$^4$ in panel C (
  $\gamma = 4.1$ pN nm$^{-1}$,
  in both panels).
  The solid, dashed and dotted curves use $N = 26,$ $N=60$ and $N = 92$ particles, and collapse
  onto a single curve when multiplied by $r_0.$
}
\label{fig:stretching}
\end{figure}
Finally, we discuss simulation data for stretching.
Consider the stretching energy of a cylindrical bilayer:
\begin{equation}
\label{eq:stretch}
\KA \frac{(A-A_0)^2}{A_0},
\end{equation}
where $A = 2\pi r L,$  $r$ is the 
midplane radius, $L$ is the cylinder length (in the $z$-direction) and $A_0$ is the area at rest.
The stretching modulus $\KA$ is for a single monolayer and 
\cref{eq:stretch} accounts for the energy of the inner
and outer monolayer leaflets of the cylinder.
Manipulation experiments give $\KA$ in the range $30$ -- $40$ $\KBT$ nm$^{-2}$
\cite{Nagle17,Nagle17-2}. 

To measure a stretching modulus, we form the circular cross
section of cylinder of radius $r$ (Figure \ref{fig:stretching}A).
The equilibrium shape is nearly circular (so long as there is a consistent number of particles in
each leaflet) and the shape obtains an equilibrium radius $r_0$ once compression and attraction are in balance. 
We use a harmonic bond to move $r$ away from equilibrium and record the change in energy
(Figure \ref{fig:stretching}BC). 

The three curves in Figure \ref{fig:stretching}B collapse onto a single curve when multplied by $r_0.$ 
Fitting to $cr_0(r-r_0)^2$ and comparing with \cref{eq:stretch}
yields $\KA =   33.4$ $\KBT$ nm$^{-2}$, $35.3$ $\KBT$ nm$^{-2}$ and $35.9$ $\KBT$ nm$^{-2}$ for the three radii
respectively.
The proximity of these three values suggests that HADF possesses a stretching modulus independent of particle number. 
Moreover, the attraction 
$\gamma = 4.1$ pN nm$^{-1}$
and repulsion
parameters $c_0 = 0.5$ pN nm$^4$ yield a consistent and physically realistic stretching modulus,
around $\KA =   35$ $\KBT$ nm$^{-2}.$
As an illustration, the curves in Figure \ref{fig:stretching}C are for the same
tension parameter and half the repulsion strength. 
There is an overall reduction in the equilibrium radii with the decreased repulsion,
and an increase in the stretching moduli $\KA$
to $40.310$ $\KBT$ nm$^{-2}$, $40.083$ $\KBT$ nm$^{-2}$ and $39.393$  $\KBT$ nm$^{-2}$
for the three curves respectively. 

\begin{figure}
\begin{center}
\includegraphics[width=\textwidth]{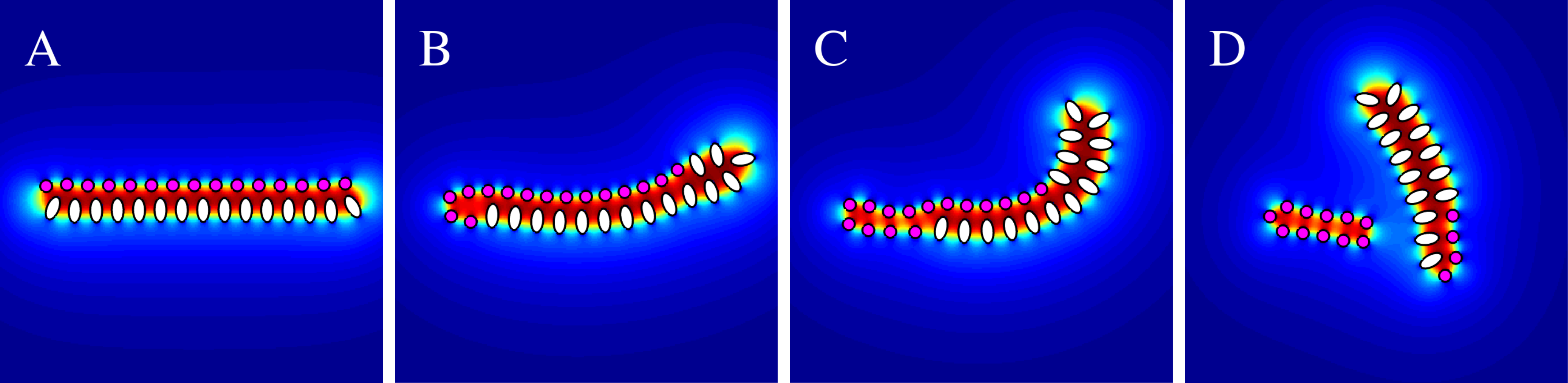}
\end{center}
\caption{Panels A--D show the spontaneous sorting and eventual fission
  of a planar bilayer mixture of small circular particles (magenta) and elliptical particles (white).
  The monolayer consisting of smaller, circular particles spontaneously migrates
  from the bilayer edge, and eventually breaks off forming forming its own
  bilayer. The bending in panels B and C suggests that the smaller, circular
particles have a spontaneous curvature more negative than for the
larger, elliptical particles. 
The parameters are 
$\gamma = 4.1$ pN nm$^{-1}$, 
$\rho   = 2.5$ nm and
$c_0    = 0.5$ pN nm$^{-4}$ 
The ellipses have $p = 6,$ $a_i = 1.25$ nm and $b_i = 0.6225$ nm
while the circles have $p = 2,$ $a_i = b_i = 0.6225$ nm.
 }
\label{fig:mixing}
\end{figure}

HADF yields physically realistic continuum-like bilayer morphologies
and these particle configurations possess elastic properties of lipid bilayer.
The HADF can also handle topological changes and mixtures in a straightforward manner.
Figure \ref{fig:mixing} illustrates the gradient descent dynamics of a lipid mixture
between small, circular and large, elliptical particles.
Under hydrophobic attraction and 
excluded volume repulsion, the particle mixture
segregates into two 
bilayers of more uniform composition. 
Diffusive interface and level-set 
approaches have dealt with the problem of mixtures by defining transport equations
for each lipid species density 
\cite{LowengrubRatzVoigt2009_PRE,MikuckiZhou17, GERA201756}. 
 
Hemifusion is one of the key intermediates
of membrane fusion involving a
Y-shaped junction between three bilayers 
\cite{Chlanda_Nature16}(see \cref{figure4}, Panel C). 
Pioneering work by Promilsow, K. and coworkers  \cite{DaiPromislow2013_RSPa,DaiPromislow2015_SIAMJMathAnal}
has lead to functionalized Cahn-Hilliard, 
diffusive interface energies that exhibit freestanding 
elastic phases, including the Y-shaped junction 
\cite{KraitzmanPromislow2018_SIAMJMathAnal,FengGuanLowengrub2018_JSciComput}. 
It is still unclear whether the HADF formulation of  
the present work is more or less efficient than a 
functionalized Cahn-Hilliard approach for capturing the 
granular energetic details of fusion \cite{Ryham16}.

\section{Hydrodynamics of amphiphilic particles in a viscous solvent}
\label{sec:dynamics}
To define particle velocities, 
we assume that the amphiphilic particles are immersed in an incompressible viscous fluid in the Stokes flow regime.
Then all the particles interact with each other through both hydrophobic forces and Stokesian hydrodynamic interactions. 
The two-dimensional particles $P_i$ have the translational and angular velocities 
\begin{equation}
\label{eq:mob_vec}
\frac{d\mathbf{a}_i}{dt} = \mathbf{v}_i,
\quad \frac{d\theta_i}{dt} = \omega_i,
\end{equation}
$i = 1,\dots,N.$
For the amphiphilic particles in a solvent, the forces 
${\bf F}_i,$ ${\bf F}_i^{\mathrm{rep}}$ 
and torques $\tau_i,$ $\tau_i^{\mathrm{rep}}$ are calculated from \cref{eq:force_torque} and \cref{eq:scalar-torque}, respectively,
The velocities $\mathbf{v}_i$ and $\omega_i$ are coupled together through the fluid velocity $\mathbf{u}$ and pressure $p$ satisfying
\begin{equation}
\label{eq:stokes}
\begin{aligned}
-\mu \Delta {\bf u} + \nabla p &= 0, \quad \\
\nabla \cdot {\bf u} &= 0,  \quad \text{in}\ \Omega,\\
{\bf u} &\to 0 \quad \text{as}\ |{\bf x}|\to \infty,\\
\mathbf{u}(\mathbf{x}) &= \mathbf{v}_i + \omega_i(\mathbf{x} - \mathbf{a}_i)^\perp,\quad \mathbf{x} \in \partial P_i,
\end{aligned}
\end{equation}
with fluid viscosity $\mu$ and subject to the stress balance conditions 
\begin{equation}
\label{eq:mobility_bc}
  \int_{\partial P_i} \mathbf{S} \cdot \mathbf{n} \,dS =  {\bf F}_i  +  {\bf F}_i^{\mathrm{rep}},\;\;\;
     \int_{\partial P_i} (\mathbf{x} - \mathbf{a}_i)\times \mathbf{S} \cdot \mathbf{n} \,dS = \tau_i + \tau_i^{\mathrm{rep}}.
\end{equation}
From  \cref{eq:net} these particle forces and torques also satisfy the force-free and torque-free conditions,
guaranteeing the existence of an integral solution for the many-body mobility problem.
The evolution equations (\ref{eq:mob_vec}--\ref{eq:mobility_bc})
satisfy the dissipation relation \cite{LEAL}
\begin{equation}
\label{eq:energy_dissipation}
\frac{d}{dt} \Phi_{\mathrm{Total}}
+ \int_{\mathbb{R}^n} \tfrac{1}{2}\mu|\nabla \mathbf{u} + \nabla \mathbf{u}^T|^2 \,dx = 0.
\end{equation}

\begin{figure}[t!]
\begin{center}
\includegraphics[width=0.99\textwidth]{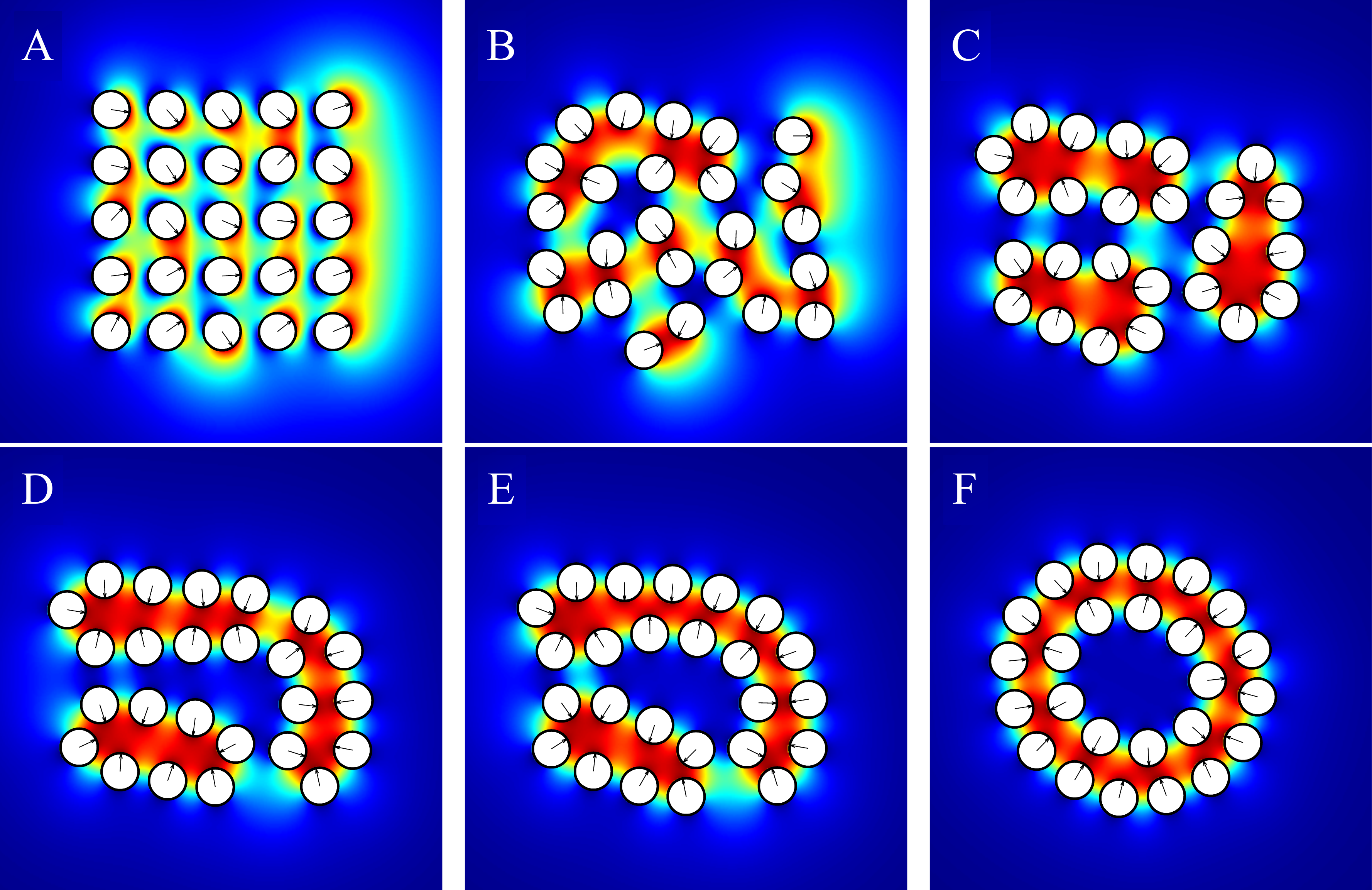}
\end{center}
  \caption{
  For large system simulations, particles self assemble into bilayer-like shapes and then eventually a cylindrical bilayer.
  Panels A--F are for $0,$  $3.8,$ $19,$ $76$, $114$ and $228$ ns, respectively;
  step size $\Delta T = 1.0[T]\approx0.38$ ns.
  The configuration in panel F evolves very slowly
  to one that is slightly more circular.
}
    \label{figure3}
\end{figure}

In two dimensions,
the kernels of single and double layer potentials for solving
the Stokes equation are the stokeslet and stresslet
\begin{equation}
\begin{aligned}
 G_{i,j}^{\mathrm{Stokeslet}}({\bf x},{\bf y}) & = \frac{1}{4\pi\mu}\bigg[\log|{\bf x}-{\bf y}|\delta_{ij}+\frac{(x_i-y_i)(x_j-y_j)}{|{\bf x}-{\bf y}|^2}\bigg], \\
 {\bf T}^{\mathrm{Stresslet}}_{i,j,k}({\bf x},{\bf y})&=-\frac{1}{\pi\mu}\frac{(x_i-y_i)(x_j-y_j)(x_k-y_k)}{|{\bf x}-{\bf y}|^4},
\end{aligned}
\end{equation}
respectively, with $i,j,k = 1,2.$
For a velocity surface density $\boldsymbol\mu,$ the stresslet satisfies the jump across the boundary 
\begin{equation}\label{eq:stokesjump}
\lim_{\mathbf{z}\rightarrow \mathbf{x}^{\pm}} \mathbf{f}_{i,\pm}(\mathbf{z})=
\mp\frac{1}{2}\boldsymbol\mu_i(\mathbf{x}) + \mbox{p.v.}\int_{\partial P_i}
{\bf T}^{\mathrm{Stresslet}}_{i,j,k}({\bf x},{\bf y})\boldsymbol\mu_j(\mathbf{y})ds_{\mathbf{y}},
\end{equation}
where $\mathbf{f}_i$ denotes the surface traction of on particle $P_i$.
Following \cite{manasthesis}, one views the external force and torque due
to hydrophobic attraction as an incident field with stress \eqref{eq:stress}.
The scattered field is then the net force and torque due to fluid mobility. 
If we split densities into 
$\boldsymbol\sigma^{\mathrm{inc}}
({\bf x})=
\{\boldsymbol\sigma^{\mathrm{inc}}_{1},\dots,\boldsymbol\sigma^{\mathrm{inc}}_{N} \}$ and
$\boldsymbol\mu({\bf x})=\{\boldsymbol\mu_{1},\dots,\boldsymbol\mu_{N} \}$, then
the particle dynamics ~\cref{eq:mobilityu} can be obtained by evaluating a single layer potential for corresponding densities $\boldsymbol\sigma$ and $\boldsymbol\mu$.
%

\begin{equation}
\begin{aligned}
\label{eq:mobilityu}
{\bf u}({\bf x}) &= \sum_{j=1}^N\int_{\partial P_j}G_{i,j}^{\mathrm{Stokeslet}}({\bf x},{\bf y})[\boldsymbol\sigma_{j}^{\mathrm{inc}}+\boldsymbol\mu_{j}]({\bf y})ds_{\mathbf{y}}\\
&= {\bf v}_i+\omega_i({\bf x}-{\bf a}_i)^\perp \quad \forall {\bf x}\in\partial P_i,
\end{aligned}
\end{equation} 
where
\begin{equation}
\begin{split}
{\bf v}_i &=\frac{d {\bf a}_i}{dt}= \frac{1}{|\partial P_i|} \int_{\partial P_i} {\bf u}({\bf y})ds_{\mathbf{y}}, \quad |\partial P_i|=\int_{\partial P_i}ds_{\mathbf{y}},\\
\omega_i=&\tau^{-1}\int_{\partial P_i}({\bf y}-{\bf a}_i)\times{\bf u}({\bf y})ds_{\mathbf{y}}, \quad \tau=\int_{\partial P_i} |{\bf y}-{\bf a}_i|^2 ds_{\mathbf{y}}.
\end{split}
\end{equation}
%

For the time-marching scheme, we solve the mobility problem for the particle translation  and rotation  velocities. We then update the particle centers and orientations using a forward Euler scheme. \Cref{alg:qbx} provides 
the time-marching details.



Non-dimensionalizing  \cref{eq:main} 
with characteristic length 1.25 nm, fluid viscosity $\mu = 1$ cP and interfacial tension $\gamma = 4.1$ pN nm$^{-1}$ 
gives the characteristic time 
$[T] =\mu a/\gamma \approx 3.82 \times 10^{-10} s.$
As an illustration, the evolution in \cref{figure1} is for 100 time steps
with time step size $\Delta T = 1.0[T].$
The time for self assembly of a few particles 
from an initially random configuration 
is thus on the order of a nanosecond. This is consistent with 
times scales for lipid rearrangements in  MD simulation \cite{Contreras10}.
Supplementary Materials Movie 1 shows the self-assembly process for three particles. 

Bilayer configurations form when we increase the number particles in the simulation. 
\cref{figure3}A has 25 Janus-type particles placed on a square grid.
The initial orientations $\theta_i$ are normally distributed about $\theta = 0$.
Within ten time steps (\cref{figure3}B), the particles rapidly rotate
to pair their hydrophobic interfaces with that of neighboring particles.  
Pairings continue to merge forming groups of eight or nine particles (\cref{figure3}C).
These groups stack together to form an arched bilayer
shape resembling the cross-section of a stomatocyte (\cref{figure3}E). \cref{figure3}F clearly shows an inner and outer monolayer  configuration a cylindrical bilayer. 

Supplementary Materials Movie 2 illustrates the self-assembly process for \cref{figure3}A--F.
As part of computational complexity test, we have calculated particle dynamics for larger systems and Supplementary Materials Movie 3 shows the results for 100 particles.

\begin{figure}
\begin{center}
\begin{tabular}{c c}
\includegraphics[width=0.35\textwidth]{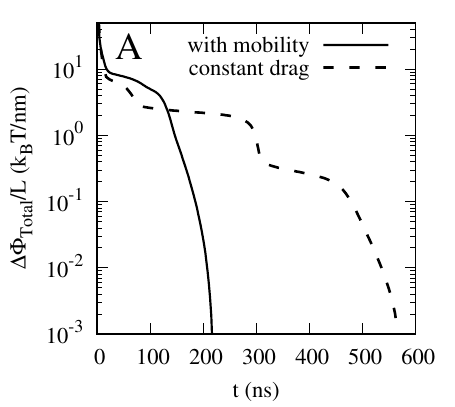}
\includegraphics[scale=0.64, trim=0 -31 0 0]{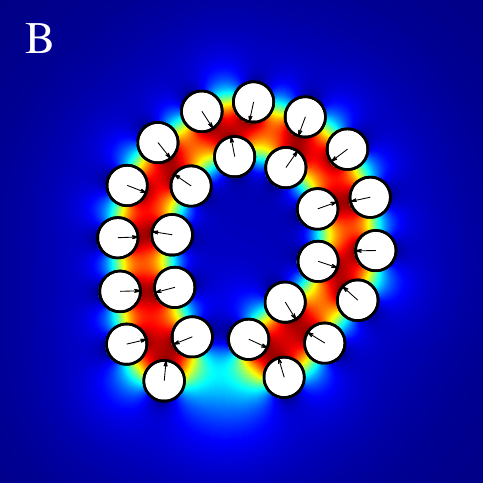} \quad
\includegraphics[scale=0.50, trim=0 -38.5 0 0]{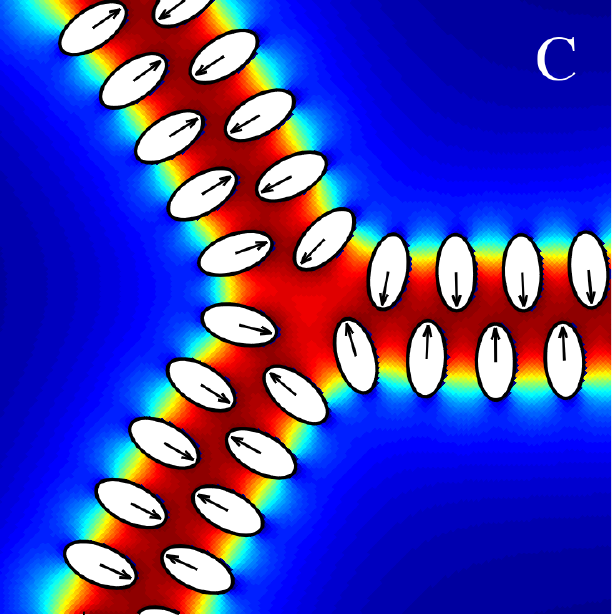} 
\end{tabular}
\end{center}
\caption{\label{figure4}
Panel A plots the change in energy $\Phi$ per length (solid curve) for the evolution \cref{figure3}A--F. 
The initial rapid decrease over 0 ns $< t <$ 10 ns corresponds to \cref{figure3}A--B.
The moderate decrease over 10 ns $< t <$ 150 ns is for \cref{figure3}C--E. 
Energy decreases very slowly over 
150 ns $< t <$ 600 ns, where
the bilayer evolves from a somewhat elliptical cross-section (\cref{figure3}F) to one that is almost perfectly circular. 
The dashed curve gives corresponding  values for a constant drag coefficient evolution: its
end-state is in panel B. 
In panel C, the free ends of three bilayers have merged into an equilibrium, Y-shaped junction. 
}
\end{figure}

Is it possible to replace the 
detailed hydrodynamic  interaction with one that uses 
a constant coefficient drag coefficient law?
The latter  also  
exhibits particle self-assembly
and avoids the computational cost 
of solving an additional mobility problem. 
Moreover,  the constant
coefficient case drag laws can closely replicate the  hydrodynamic interaction case when 
particles are dispersed (see Supplementary Material, Section S2.3).
Nevertheless, numerical experiments show that 
the choice of dissipative mechanism is consequential 
to the time course. For example, 
\cref{figure3}F and \cref{figure4}B compare the two 
different end-states resulting from a identical initial configurations for a viscous fluid and for a constant drag law, respectively. The difference lies in the number of particles contained in the inner leaflet and this  has determined whether or not the bilayer closes.


\section{Boundary Integral Equation Formulation}
\label{sec:bie}

In this section, we present a second kind integral equation formulation
for the exterior Dirichlet problem
\cref{SL}.
The domain $\Omega$ is the exterior domain, meaning that its complement $\Omega^c$ (the collection of particles) is compact.

There are a number of numerical methods for solving the exterior problem \cref{SL}. These 
include finite difference methods, finite element methods, and boundary integral equation (BIE) methods. 
The BIE methods are perhaps the most suitable
since they represent the solution via layer
potentials with an unknown density only on the boundary. 
This reduces the dimension of the problem by one and leads to a much smaller linear system. 
Another advantage is that the integral representation automatically satisfies the
governing partial differential equation and the
boundary condition at infinity. Thus, there is no need
to truncate the computational domain and impose artificial boundary conditions,
as would be the needed with the finite element and finite difference approaches. Finally, 
when combined with high-order quadratures and fast algorithms such as the fast
multipole methods \cite{fmm,Greengard02}, the BIE formulation leads to a high-order numerical algorithm with optimal computational complexity.

Before describing the method, we first consider whether  the far-field condition in \cref{SL} is sufficient to 
determine a unique solution. 
As mentioned in Section~\ref{sec:intro},
the functional $I[\cdot]$ has a unique minimizer in $\mathcal{A}$.
The minimizer $u$ satisfies 
$0 \leq u(x) \leq 1$ for all $x \in \Omega$. 
To see why these bounds holds, 
consider a truncated version  $\tilde u = \max\{0,\min\{1, u\}\}$ of $u$.
Because $0 \leq f \leq 1$ on $\Sigma$, $\tilde u\in\mathcal A$.
Lastly, $I[\tilde u] \leq I[u]$ by inspection, and
$I[u] \leq I[\tilde u]$ by minimality of $u$. This implies that $I[\tilde{u}] = I[u],$
and since $u$ is the unique minimizer, we have  $u = \tilde u$. 

To obtain a far-field decay condition, select a sufficiently large distance $D > 0$ from
$\Omega^c$. Let $x$ be such that $d(x, \Omega^c) > D$. By a change of
coordinates, we may assume that $x$ is the origin $(0, \dots, 0)$ and that $\Omega^c$ lies in the set
$\{ (x_1,\dots, x_n): x_n \leq - D \}$.  In this coordinate system, consider the function
\[
  v_{\epsilon}(x_1,\dots,x_n) = \frac{\cosh(x_n / \rho)}{\cosh(D / \rho)} + \epsilon
  \sum_{i=1}^{n-1}\cosh(x_i / \rho), \quad \epsilon > 0.
\]
Then $-\rho^2 \Delta v_{\epsilon} + v_{\epsilon} = 0$ and $v_{\epsilon} \geq 1$ on $\Sigma$ (since $x_n$ is less than $-D$ there), and $v_{\epsilon} \to \infty$
as $|(x_1,\dots, x_n)| \to \infty$.
Next, consider the function $g = v_{\epsilon} - u$. We have $g \geq 0$ on $\Sigma$ and $g \geq 0$
for all sufficiently large $(x_1,\dots, x_n)$ (since $v_{\epsilon} \to \infty$ and 
$u$ is bounded everywhere between $0$ and $1$).

From the weak maximum principle~\cite[Cor.~3.2, assuming $f \in C(\Sigma)$]{gilbarg_elliptic_2015}, we have $g \geq 0$ in $\Omega$. It follows that $v_{\epsilon} \geq u$ in
$\Omega$.  Finally, letting $\epsilon \to 0$, we conclude that
\begin{equation}
   %
   %
   0 \leq u(x) \leq v_{\epsilon}(x) < 2 \exp(-D / \rho)
   \label{exp_decay}
\end{equation}
as soon as $\epsilon > 0$ is sufficiently small.

The problem \cref{SL} thus has at least one solution (vanishing at infinity), namely the variational one. But since the domain $\Omega$ is non-compact,
it is in principle possible that \cref{SL} has multiple solutions vanishing
at infinity with different rates. The following Liouville-type result shows that 
this is not the case. In fact, we get uniqueness even if we replace the zero
far-field condition with the power growth condition $u = O(|x|^p)$
as $x \to \infty$.

\begin{lemma}
\label{liouville}
The exterior problem \cref{SL} has at most one solution. 
\end{lemma}
\begin{proof}
Suppose that \cref{SL} has two solutions $u_1$ and $u_2$. 
Let $w = u_1 - u_2$ and define 
\begin{equation*}
E(r) = \int_{\Omega \cap \mathbf{B}_r} \tfrac{1}{2}w^2 \,dx,\qquad E: [0,\infty) \to [0,\infty),
\end{equation*}
where $B_r$ is the ball of radius $r$ centered at the origin. 
Select $r_0$ positive and sufficiently large so that $\Omega^c \subset B_{r_0}$.

The function $E(r)$ is infinitely differentiable on $(r_0,\infty)$ 
since any solution of \cref{SL} is smooth in $\Omega.$ 
Using $w = 0$ on $\Sigma$ and \cref{SL}, it follows that 
\begin{equation}
E''(r) = \frac{n-1}{r}E'(r) + \frac{1}{\rho^2} E(r) + \int_{\Omega \cap B_r} |\nabla w|^2\,dx,\quad r_0 < r < \infty.
\end{equation}
Since $E$ is nondecreasing, by definition, we get 
\begin{equation*}
E''(r) \geq \rho^{-2} E(r),\quad r_0 < r < \infty. 
\end{equation*}

Let $F(r) = \tfrac{1}{2}E(r_0) \exp((r-r_0)/\rho).$ We claim that 
\begin{equation}
\label{growth}
F(r) \leq E(r), \quad r_0 \leq r < \infty.
\end{equation}
To form the comparison argument, suppose to the contrary that 
$E(r)$ is not everywhere greater than or equal to $F(r)$. 
Then $E(r_0) > 0$ and there is $r_* > r_0$ 
with $E'(r_*) \leq F'(r_*)$ and $F(r) < E(r)$ for $r_0 \leq r < r_*$.
But then,
\begin{equation*}
E'(r_*) \leq F'(r_*) 
= \int_{r_0}^{r_*} \rho^{-2} F(r) \,dr 
< \int_{r_0}^{r_*} \rho^{-2} E(r) \,dr
\leq  E'(r_*) - E'(r_0).
\end{equation*}
These inequalities are in contradiction since $E'(r_0)$ is nonnegative.

Thus, if $u_1$ and $u_2$ are two solutions 
of \cref{SL}, then according to \cref{growth} there are 
one of two possibilities: either $E(r_0) = 0$  
or one of the solutions has exponential growth.
The vanishing condition $u \to 0$ as $x \to \infty$ rules out the latter case
and so $E(r_0)$ must be zero. Since $r_0$ was arbitrary, $u_1$ and $u_2$ are identical. 
\end{proof}
 
In two dimensions, the equation $-\rho^2 \Delta u + u=0$
has the free-space Green's function (also called fundamental solution) 
\begin{equation}
        G(x, y) = \frac1{2\pi}K_0(|x-y|/\rho), \quad x,y \in \mathbb{R}^2,
        \label{Gfunc}
\end{equation}
where $K_0$ is the zeroth order modified Bessel function
of the first kind~\cite{nisthandbook}.
For a Lipschitz domain $\Omega$ in $\mathbb{R}^2$ with boundary $\Sigma$, the space
$L^2(\Sigma)$ denotes all square integrable functions on $\Sigma$. 
Given a function $\sigma\in L^2(\Sigma)$, we define the single layer potential by the formula
\begin{equation}
\mathcal{S}[\sigma](x) = \int_{\Sigma}G(x,y)\sigma(y)ds_y,
\label{eq:single-layer-def}
\end{equation}
and the double layer potential by the formula
\begin{equation}
\mathcal{D}[\sigma](x)=\int_{\Sigma} \frac{\partial G(x,y)}{\partial \nu(y)} \sigma(y)ds_y,
\label{eq:double-layer-def}
\end{equation}
where $\nu(y)$ is the unit outward normal vector
with respect to $\Omega^c$.
It is well-known from classical potential 
theory~\cite{kress2014} that
the single layer potential is continuous and the double layer potential exhibits a jump across the boundary.
To be more precise, when $z$ approaches a point  $x\in\Sigma$ nontangentially, the limits of
$\mathcal{S}[\sigma]$ and $\mathcal{D}[\sigma]$ exist and are given by the following formulas:
\begin{equation}\label{sjump}
\lim_{z\rightarrow x^{\pm}} \mathcal{S}[\sigma](z) =
S[\sigma](x) = \int_{\Sigma}G(x,y)\sigma(y)ds_y,
\end{equation}
and
\begin{equation}\label{eq:djump}
\lim_{z\rightarrow x^{\pm}} \mathcal{D}[\sigma](z)=
(\pm\frac{1}{2}I+ D)[\sigma](x)=
\pm\frac{1}{2}\sigma(x) + \mbox{p.v.}\int_{\Sigma}
\frac{\partial G(x,y)}{\partial \nu(y)} \sigma(y)ds_y,
\end{equation}
for almost every point $x\in \Sigma$. Here $z\rightarrow x^{\pm}$
implies that $z$ approaches $x$ from the exterior$(+)$ or the interior$(-)$
of $\Omega^c$, respectively. It is also well-known
that both the single layer operator
$S: L^2(\Sigma)\rightarrow L^2(\Sigma)$ and
the double layer operator
$D:  L^2(\Sigma)\rightarrow L^2(\Sigma)$ are compact when the boundary $\Sigma$ is $C^1$.

We will represent the solution to \cref{SL}
with the double layer potential representation:
\begin{equation}
        u(x) = \mathcal{D}[\sigma](x).
        \label{intrep}
\end{equation}
The jump relation of the double layer
potential~\cref{eq:djump} leads
to the following boundary integral equation
on the unknown density $\sigma$:
\begin{equation}
        \frac12\sigma(x)+D[\sigma](x) =f(x), \quad x \in \Sigma.
        \label{bie}
\end{equation}
\begin{algorithm}[ht]
\caption{Particle Updates by Exterior Screened Laplace BVP}
\label{alg:qbx}
\begin{algorithmic}[1]
\STATE{Set the particle centers $\mathbf{a}_i \in \mathbb{R}^2$ and orientations $\theta_i \in \mathbb{R}$,  boundary condition $f(x)$ and time step size $\Delta t$}
\STATE{Determine the discretization on the boundary $\Sigma$ and construct the double layer potential $\mathcal{D}[\sigma](x)$.}
\FOR{$t = t_0 : t_{end}$}
\STATE{Use GMRES iterative method to solve the unknown density $\sigma $ in~\cref{bie}.}
\STATE{Use the solved $\sigma$ to obtain the screend Laplace equation solution $u$ 
and calculate $\nabla u$.}
\STATE{Calculate inter-molecular forcesr~\cref{eq:force_torque} and \cref{repul}.}
\STATE{Use GMRES iterative method to solve the unknown density $\boldsymbol\mu$ in~\cref{eq:stokesjump}. }
\STATE{Solve mobility problem
 and update particle velocities ${\bf v}_i$ and $\omega_i$ in \cref{eq:mobilityu}.}
\STATE{Update particle center positions ${\bf a}_i$ and orientations $\theta_i$.}
\STATE{Update the marching time $t = t_0 + \Delta t$.}
\ENDFOR
\RETURN $T$
\end{algorithmic}
\end{algorithm}

\begin{theorem}
Suppose that $\rho$ is any positive real number. Then for any $f\in L^2(\Sigma)$, the 
second kind integral equation \cref{bie}
is uniquely solvable.
\end{theorem}
\begin{proof}
By the Fredholm alternative~\cite{kress2014}, we only
need to show that the only solution to
the homogeneous equation
\begin{equation}
        \frac12\sigma(x)+D[\sigma](x) =0.
        \label{biehomo}
\end{equation}
is $\sigma\equiv 0$.

Consider the function $u(x)$ defined by 
the formula \cref{intrep}. It is clear that
$u$ satisfies the equation $-\rho^2 \Delta u + u = 0$ in both
the exterior domain $\Omega$ and the interior
domain $\Omega^{c},$ and vanishes at infinity. By the uniqueness
of the exterior Dirichlet problem (Lemma \ref{liouville}), we have $u\equiv 0$ in $\Omega$. Hence,
\begin{equation}
\lim_{z\rightarrow x^{+}}\frac{\partial u(z)}{\partial \nu}=0, \quad x\in\Sigma.
\end{equation}
Since the normal derivative of the double layer
potential is continuous across the 
boundary~\cite[gen. of Thm.~6.18]{kress2014}, we have
\begin{equation}
\lim_{z\rightarrow x^{-}}\frac{\partial u(z)}{\partial \nu}=0, \quad x\in\Sigma.
\end{equation}
Hence, $u$ in the interior domain $\Omega^{c}$
is the solution to  
 the interior Neumann problem
\begin{equation}\label{intneum}
-\rho^2 \Delta u + u=0, \quad u\in \Omega^{c},\quad 
\frac{\partial u}{\partial \nu}=0, \quad x\in \Sigma.
\end{equation}
Applying Green's first identity, we obtain
\begin{equation}
\int_{\Omega^{c}} \rho^2 |\nabla u|^2 + u^2 dx = 0.
\end{equation}
Thus we have $u\equiv 0$ in $\Omega^{c}$ as well. 
The jump relation of the double layer potential \cref{eq:djump} leads to
\begin{equation}
\sigma(x)=\lim_{z\rightarrow x^{+}} u(z)-
\lim_{z\rightarrow x^{-}} u(z)=0, \quad x\in \Sigma,
\end{equation}
which completes the proof.
\end{proof}

\begin{remark}
As pointed out earlier, the screened Laplace equation 
can be viewed as the Helmholtz equation
$\Delta u+k^2 u=0$ with pure imaginary $k$.
When $k$ is an arbitrary complex number, the so-called Brakhage-Werner representation~\cite{Brakhage1965} (also called the Burton-Miller
representation~\cite{burton1971prs} in acoustics) represents the solution to the Helmholtz equation via a linear combination of single and double layer potentials
\begin{equation}
        u(x) = i\mathcal{S}[\sigma](x) + \mathcal{D}[\sigma](x).
        \label{eq:bwform}
\end{equation}
It has been shown that the  representation \cref{eq:bwform} leads to a uniquely solvable second kind integral equation for any value of 
$k \in \mathbb{C}$~\cite{nedelec2001}. Due to 
the exponential decay of the solution to our
exterior problem \cref{SL}, we are able to use the double layer potential alone to represent
its solution and still achieve existence and uniqueness of the associated boundary integral equation \cref{bie}.
\end{remark}

We would like to point out that when $\rho$ is very large, the formulation~\cref{eq:bwform} may 
lead to a better conditioned linear system
than \cref{intrep}. There are other second kind integral equation formulations for this problem. For example, one may replace the single layer potential by a collection
of point sources inside each particle, where the strength
of the point source may be unknown or equal to the average value of the unknown density
function on the boundary of each particle. We refer
the readers to~\cite{greenbaum1993jcp,greengard1996jcp,jiang2017siammms}
for details.


\section{High-order quadrature and fast algorithms}
\label{sec:numerical_methods}

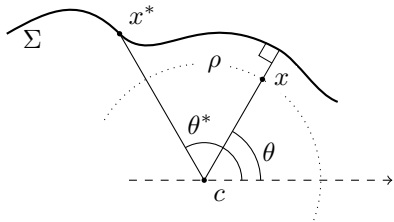
\begin{figure}
  \centering
  \begin{tikzpicture}
    \coordinate (c) at (0,0) ;
    \path (c) ++(45+75:2.25) coordinate (s);
    \path (c) ++(45+15:1.55) coordinate (t);
    \path (c) ++(45+15:2) coordinate (t-to-curve);

    \path (s) ++(-1.5,0) coordinate (curve-before);
    \path (t) ++(1,-0.3) coordinate (curve-after);
    \draw [thick]
      (curve-before)
      ..controls +(30:0.7) and +(180-45:0.7) ..
      (s)
      node [pos=0.2,anchor=north] {$\Sigma$}
      ..controls +(-45:0.7) and +(45+90+15:1.25) ..
      (t-to-curve)
      ..controls +(45+270+15:0.3) and +(160:0.3) ..
      (curve-after) ;

    \fill (c) circle (1pt);
    \fill (t) circle (1pt);
    \fill (s) circle (1pt);

    \draw [dotted] ++(-20:1.55) arc (-20:45+105:1.55);

    \node at (85:1.55) [fill=white] {$\rho$};
    \draw (c) -- (t) ;
    \draw (c) -- (s) ;
    \node at (c) [anchor=north west] {$c$};
    \node at (s) [anchor=south west] {$x^*$};
    \node at (t) [anchor=west,fill=white,xshift=0.5mm,inner sep=1mm] {$x$};
    \draw [dashed,->] (c) ++(-1, 0) -- ++(3.5,0);

    \draw (c) ++(0.5,0) arc (0:45+75:0.5);
    \path (c) ++(45*0.8+75*0.8:0.5)
      node [anchor=south] {$\theta^*$};

    \draw (c) ++(0.75,0) arc (0:45+15:0.75) ;
    \path (c) ++(45/2+15/2:0.75)
      node [anchor=west] {$\theta$};

    \draw [very thin] (t) -- (t-to-curve);
    \draw [very thin]
      let
        \p1 = ($ 0.1*(t-to-curve) - 0.1*(c) $),
        \p2 = (-\y1,\x1)
      in
      ($(t-to-curve)!0.1!(c)$) -- ++(\p2) -- ++(\p1) ;

  \end{tikzpicture}
  \caption{%
    Geometric situation of a single QBX expansion
    with sources along the collection of curves $\Sigma$, as
    used in~\cref{eq:local-exp} and~\cref{eq:graf-coefficient}.
    Note that the target point $x$ will reside on $\Sigma$ for the
    computation of the on-surface value of the layer potential.
  }%
  \label{fig:qbx-geometry}
\end{figure}
For the accurate and rapid evaluation of the layer potentials
occurring in the previous section, we make use of `Quadrature by
Expansion', or QBX for short~\cite{Andreas}. Here we briefly review the QBX scheme for 
the general Helmholtz kernel and note that the Yukawa kernel  \cref{Gfunc} is simply a special case of the Helmholtz kernel with pure imaginary wave number $k$. To do so, we cover
a neighborhood of the source curve $\Sigma$ with locally valid
(`local') expansions of the potential emanating from the entire source
curve $\Sigma$. For a collection of on-surface target points $(x_i)$,
expansion centers are chosen as $c_i=x_i+\nu \eta_{x_i}$, where $\eta_x$ is
a scaling factor connected to the local quadrature resolution.
See~\cite{gigaqbx3d} for details of the determination of $\eta_x$.
Then, for target points $x\in B(c_i,\eta_{x_i})$, the layer potential
may be evaluated as
\begin{equation}
 \phi(x) = \sum_{l=-\infty}^\infty \alpha_l J_l(k \rho) e^{-i l \theta}
 \label{eq:local-exp}
\end{equation}
where $(\rho,\theta)$ denote the polar coordinates of the target point $x$
with respect to the expansion center $c$, and $J_l$ is the Bessel
function of order $l$ (see Fig.~\ref{fig:qbx-geometry}).  For the
single layer potential $S\sigma$, the coefficients $\alpha_l$ in the
expansion (\ref{eq:local-exp}) can be computed analytically:
\begin{equation}
  \alpha_l = \frac{i}{4} \,
  \int_\Sigma H^{(1)}_l(k|x^*-c|) e^{i l \theta^*} \sigma(x^*) d x^*,
  \qquad
  (l = -p, -p+1, \ldots, p)
  \label{eq:graf-coefficient}
\end{equation}
where $(|x^*-c|,\theta^*)$ denote the polar coordinates of the point $x^*$ with respect
to $c$. These (now non-singular) integrals for the coefficients
$\alpha_l$ are then computed by conventional high-order numerical
quadrature.
These formulas follow immediately from
Graf's addition theorem~\cite[(10.23.7)]{nisthandbook},
\begin{equation}
  H^{(1)}_0(k|x-x^*|) = \sum_{l=-\infty}^\infty H^{(1)}_l(k|x^*-c|) e^{i l \theta^*}
  J_l(k |x-c|) e^{-i l \theta},
  \label{eq:graf-addition}
\end{equation}
This identity applies directly to the Yukawa potentials under consideration
here, based on the fact that $K_0(z) =(i\pi/2) H_0^{(1)}(iz)$,
cf.~\cite[(10.27.8)]{nisthandbook}.
Separation-of-variables results similar to Graf's addition theorem
hold for Laplace potentials, allowing us to proceed analogously
in that case~\cite{fmm}. The QBX procedure described above employs two
means of approximation: the truncation of the series expansion, and
the computation of the coefficients by numerical quadrature.
We give an error result for QBX that accounts for both aspects.  For
the following result, we consider the case of the double layer and
assume $c=0$ without loss of generality.

\begin{theorem}[{QBX truncation and quadrature
errors,~\cite[Thm.~2.5 and (4.6)]{epstein:qbx-error-est}}]
  \label{thm:qbx-accuracy}
  Suppose that $\Gamma$ is a smooth, bounded curve
  embedded in $\mathbb R^2$, such that $B_{\eta_x}(0)\cap\Gamma=\emptyset$, but
  $\rho e^{i\theta}\in \overline B_{\eta_x}(0)\cap\Gamma$.
  Assume the geometry $\Gamma$ is discretized using $q$ point composite
  Gauss-Legendre panels of uniform length $h$, with a total of $n$ points.

  For $k\in[0,\infty)$, $N$ a 
  positive integer, and $\beta>0$, there are constants
  $C'_{N,\beta,\Gamma}(k)$
  and  $C''_{q,\Gamma}(k)$,
  so that if $\sigma\in C^{N,\beta}(\Gamma)$, then
  \begin{multline}
    \left|
      \lim_{r \to \rho^-}  \int_{\Gamma}
      \frac{\partial G(r e^{i\theta},y)}{\partial \nu(y)}    
      \sigma(y)ds(y)
      -\sum_{l=1-N}^{N-1} Q_q(\alpha_l) J_{l}(k\rho)e^{-il\theta}
    \right| \\
    \leq
    C'_{N,\beta,\Gamma}(k)\rho ^{N}\|\sigma\|_{C^{N,\beta}(\Gamma)}
    + C''_{q,\Gamma}(k) \frac{h^{2q}}{{(4r)}^{2q+1}} \| \sigma \|_{C^{2q}}.
  \end{multline}
  Here the coefficients $\{\alpha_l\}$ are given by~\cref{eq:graf-coefficient},
  and $Q_q(\alpha_l)$ denotes the approximation of the coefficient
  integral of~\cref{eq:graf-coefficient} by Gaussian quadrature with $q$ points.
\end{theorem}

The theorem makes several assumptions that may not be true of the
geometry discretization in its original form, notably the assumption
that the placed disks do not intersect with other geometry, or the
requirement that source panels supply sufficient quadrature
resolution not just for themselves, but also for adjacent panels
(which masquerades in Theorem~\ref{thm:qbx-accuracy} as the
assumption of equal panel sizes). All these issues can be remedied by
adequate refinement of the source geometry. An efficient, tree-based
algorithm is available~\cite{gigaqbx3d} to accomplish this.

To avoid quadratic scaling of the computational cost with the number
of degrees of freedom, boundary integral equation methods require some
form of acceleration, often through a variant of the Fast Multipole
method (FMM~\cite{carrier:1988:adaptive-fmm}). In the context of
QBX, it is convenient to exploit that the expansions produced
as the output of the far-field stage of the FMM are the same ones
employed by the quadrature method. However, without some
care, loss of accuracy may occur~\cite{rachh:2017:qbx-fmm}.
We make use of the `GIGAQBX' fast algorithm of~\cite{gigaqbx2d}
to obtain guaranteed accuracy at linearly scaling cost. This algorithm
modifies the conventional FMM by forcing direct computation of
interactions that may endanger the accuracy of the computed
QBX expansion, in addition to a number of modifications to retain
efficiency and linear scaling in that setting.

Another feature in GIGAQBX is the adaptive refinement activated 
when two or more source geometries get close to each other, causing near-singular evaluations of the boundary integral.
The adaptive refinement is designed to continue until the expansion 
disks get out of the region of the source geometry~\cite{gigaqbx2d}. Under the HADF,
one might expect that many levels of refinement are needed when two particles are brought to near contact. However, the steric potential~\cref{repul} also acts to prevent particles from getting too close to each other.

With the use of the short range repulsive potential, we found that the count of continuous refinements to be at most 3 to 5 levels at each time step for simulations presented in this work. 
Moreover, if the target point is geometrically on the wrong side (e.g. the interior region for the exterior problem), the GIGAQBX approximates analytic continuation of the potential across the boundary $\Gamma$,
leading to benign behavior even in degenerate cases.

Our simulation codes make use of the software package
`\textsc{Pytential}'~\cite{pytential}, which is in turn built upon
\textsc{FMMlib}~\cite{fmmlib} for some of its expansion and
translation operator infrastructure.


\section{Conclusions}
\label{sec:conclusions}

Topological transitions of a lipid bilayer membrane, such as membrane fusion and fission, involve rearrangement of lipid molecules in the bilayer. Consequently the well-known Helfrich free energy requires modification to account for lipid granularity 
to resolve the detailed lipid re-modeling during membrane fusion or fission \cite{OhtaKawasaki,DaiPromislow2015_SIAMJMathAnal,Ryham16}. 
By using a modified Helfrich free energy with van der Waals repulsion and a hydrophobic potential for lipid tail-solvent
interaction, Ryham {\it et al.}~\cite{Ryham16} calculated a least energy pathway of membrane fusion.
Building on these results, the main motivation for the work presented here is to construct a hybrid continuum lipid model at the mesoscopic scales to capture 
both the lipid granularity and the long-range interaction during 
the self-assembly of lipid molecules and fusion/fission dynamics of a lipid bilayer membrane.

Our continuum coarse-grained model for lipids focuses on the
hydrophobic interactions between lipid tails, and an SKIE
formulation of the hydrophobic stress is derived and used for
obtaining particle dynamics. We also show that the long-range hydrophobic
attraction potential is non-pairwise, and thus requires special treatment within the coarse-grained model
framework. The GIGAQBX scheme---an improved version of the QBX-FMM scheme with guaranteed accuracy---is used in the discretization, solver, and evaluation phases of the SKIE to achieve high accuracy and asymptotically optimal complexity.
Simulation results of our model show that during the self-assembly
process, coarse-grained lipid particles form structures (such as
micelles and bilayers) that may further fuse together to form a single bilayer
membrane. These results show that our approach can naturally capture the mesoscopic dynamics of
membrane fusion/fission. Furthermore,
we show that the hydrophobic interactions give rise to membrane
curvature minimization, which is an indication of the origin of bending rigidity in a bilayer membrane.

It is straightforward to apply the numerical scheme
developed in this paper to study particles of arbitrary
shape. With slight algorithmic modification, the scheme
can also accurately capture the collision dynamics that many researchers may regard as rather difficult to deal with. 

We also illustrate that the lipid hydrodynamics under HADF gives rise to macroscopic mechanical properties of a lipid bilayer membrane that are consistent with
other results in the literature. The flexibility of our hybrid approach allows us to consider a mixture of two lipid species and how spontaneous sorting (phase separation) of two
lipid species leads to membrane fission, consistent with results from phase-field simulations \cite{LowengrubRatzVoigt2009_PRE}.
Our future goal is to extend the current framework to three-dimensional lipid system. We will incorporate fluctuating hydrodynamics into the boundary
integral formulation to extract physical properties of the lipid bilayer membrane such as membrane diffusivity, bending rigidity and the surface tension.
By modification of the interfacial labels, 
HADF can account for charged lipids and study their impact on elastic
properties of bilayer \cite{Soft_matter_Dimova_Vlahovska2019}.
Since we have immersed the particles in a zero-Reynolds flow, it
is possible to study the rheological properties of micelle networks 
in large particle simulations \cite{Lutz_Bueno_Langmuir16}.
Finally, we also aim to investigate the continuum limit of our hybrid model and make comparison with functionalized Canham-Helfrich models.


\section*{Acknowledgments}
 Part of the work was performed when
S.-P.~Fu, S.~Jiang, A.~Kl{\"o}ckner, and M.~Wala
were participating the 2017 HKUST-ICERM workshop 
"Integral Equation
Methods, Fast Algorithms and Their Applications to Fluid Dynamics and
Materials Science."
The authors thank the anonymous referees,
Joshua Schrier, Jasun Gong, John Lowengrub and Jun Allard for valuable feedback. S. Jiang was supported by NSF under grants
DMS-1418918 and DMS-1720405, and by the Flatiron
Institute, a division of the Simons Foundation. A. Kl{\"o}ckner and M.~Wala were supported in part by NSF under grants DMS-1418961 and DMS-1654756.
Y.-N. Young was supported by NSF under grants DMS-1412789 and DMS-1614863.



\begin{appendix}
\renewcommand{\arraystretch}{1.3}

\begin{figure}[ht]
\begin{center}
\includegraphics[width=0.24\textwidth,  trim=0 66 0 0]{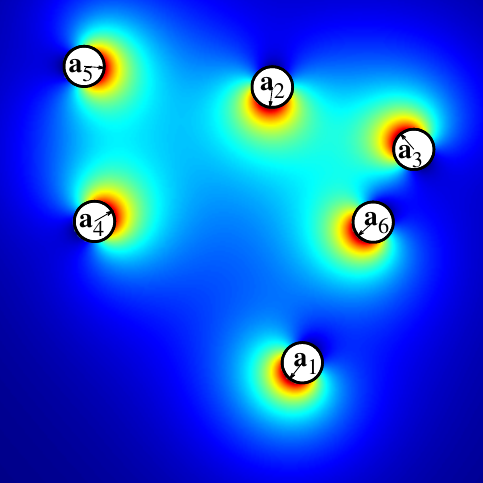}\footnotesize
\quad 
    \begin{tabular}{ |r| r| r|r| r| c|}
    \hline
\thead{$i$}   & \thead{${\bf F}_i^x$} &  \thead{${\bf F}_i^y$} & \thead{$\mathbf{f}_{i}^x$} & \thead{$\mathbf{f}_{i}^y$} & \thead{rel. diff. (\%)}  \\ \hline
	$1$ &  $0.017$  &  $0.124$ & $0.022$ & $0.177$ & $42.6\%$   \\
	$2$ &  $0.095$  & $-0.167$ & $0.122$ & $-0.262$ & $51.5\%$  \\
	$3$ & $-0.474$  & $-0.019$ & $-0.555$ & $-0.075$ & $20.9\%$   \\
	$4$ &  $0.173$  &  $0.255$ & $0.231$ & $0.291$ & $22.0\%$    \\
	$5$ &  $0.154$  & $-0.298$ & $0.208$ & $-0.338$ & $20.1\%$    \\
	$6$ &  $0.035$  &  $0.104$ & $-0.027$ & $0.207$ & $109.5\%$    \\
    \hline
    \end{tabular}
\end{center}
  \caption{The pseudo-color map shows a configuration of six,
  randomly placed particles with random orientations. 
  The table (in pN) provides values for the $x$ and $y$ components
  of force ${\bf F}_i^{\{x,y\}}$ 
  calculated from the HADF and the force ${\bf f}_i^{\{x,y\}}$ 
  calculated assume a pairwise potential \cref{pairwise_assumption}. 
  The rightmost column shows the relative difference 
$\|\mathbf{F}_i - \mathbf{f}_i\|/\|\mathbf{F}_i\|.$
  }
    \label{fig:nonpair}
\end{figure}
\renewcommand{\arraystretch}{1.0}

\section{Pairwise potentials}
\label{sec:non-pairwise}
We show that pairwise potentials do not 
closely approximate the HADF.
Consider the case of $N$ many particles in general position and orientation. Their associated pairwise potential is 
\begin{equation}
\label{pairwise_assumption}
\sum_{i=1}^N \sum_{j > i} \phi_{ij},
\end{equation}
where $\phi_{ij}$ is the functional \eqref{eq:main} evaluated
on $\Omega = \mathbb{R}^2 \setminus (P_i \cup P_j).$
Differentiating \cref{pairwise_assumption} with respect to  
position $\mathbf{a}_i$ yields the force
\begin{equation}
\label{pairwise_force}
\mathbf{f}_i =  \sum_{j\neq i} \mathbf{f}_{ij},\quad \mathbf{f}_{ij} = -\nabla_{\mathbf{a}_i} \phi_{ij}.
\end{equation}
That is, we calculate $\mathbf{f}_{ij}$ using \cref{eq:force_torque} 
for a fluid domain containing only two particles, $P_i$ and $P_j,$ 
and then sum the results for $j = 1,\dots, N,$ $j \neq i.$
Finally, let $\Phi$ be the HADF for all $N$ particles,
and calculate the hydrophobic force $\mathbf{F}_i = -\nabla_{\mathbf{a}_i} \Phi$ using \cref{eq:force_torque}
over the fluid domain that contains all particles. 

The table in Figure \ref{fig:nonpair} compares 
the non-pairwise $\mathbf{F}_i$ and pairwise $\mathbf{f}_i$ forces
for a sample particle configuration 
(Figure \ref{fig:nonpair}, pseudo-color map).
The forces show significant differences for all 
six particles (Figure \ref{fig:nonpair}, rightmost column),
suggesting that it is insufficient to use a pairwise potential to 
calculate HADF as formulated in the present work. 
We note, however, that owing to the form
of the free-space Green's function \cref{Gfunc}, the correlations between particles decays like $\exp(-D/\rho)$ 
in their distance $D.$ This makes it possible
localize interaction to tens of particles 
by setting a cut-off radius in the layer potential
evaluations. 

\end{appendix}

\bibliographystyle{siamplain}
\bibliography{references}
\resetlinenumber[1]

\headers{Supplementary Material: Amphiphilic Interaction}{S.-P. P. Fu, R. J. Ryham, A. Kl{\"o}ckner, M. Wala, S. Jiang, and Y.-N. Young}



\thispagestyle{empty}

\newpage
{\Large \bf

  \noindent Supplementary Material\\

  \noindent 
  Simulation of Multiscale Hydrophobic Lipid Dynamics via Efficient Integral Equation Methods}\\

\noindent 
Szu-Pei P. Fu$^{1,*},$ 
Rolf J. Ryham$^{1},$ 
Andreas Kl{\"o}ckner$^{2},$ 
Matt Wala$^{2},$
Shidong~Jiang$^{3},$,
Y.-N. Young$^{3},$
\\

\noindent
$^{1}$Fordham University, Department of Mathematics,  Bronx, NY, USA

\noindent
$^{2}$Department of Computer Science, University of Illinois at Urbana-Champaign, Urbana, IL 61801 USA

\noindent
$^{3}$Department of Mathematical Sciences, New Jersey Institute of Technology, Newark, NJ  07102 USA
\\

\noindent $^*$Corresponding author. Address: Fordham University, Department of Mathematics, 441 E. Fordham Rd, Bronx, NY 10458. email: \text{sfu17@fordham.edu}

\setcounter{page}{1}

\setcounter{figure}{0}
\renewcommand{\thefigure}{S\arabic{figure}}

\setcounter{equation}{0}
\renewcommand{\theequation}{S\arabic{equation}}

\setcounter{section}{0}
\renewcommand{\thesection}{S\arabic{section}}



\newpage

\sloppy
\section{Elliptical Excluded Volume Repulsion}

\begin{figure}[h]
\begin{center}
\includegraphics[width=0.4\textwidth]{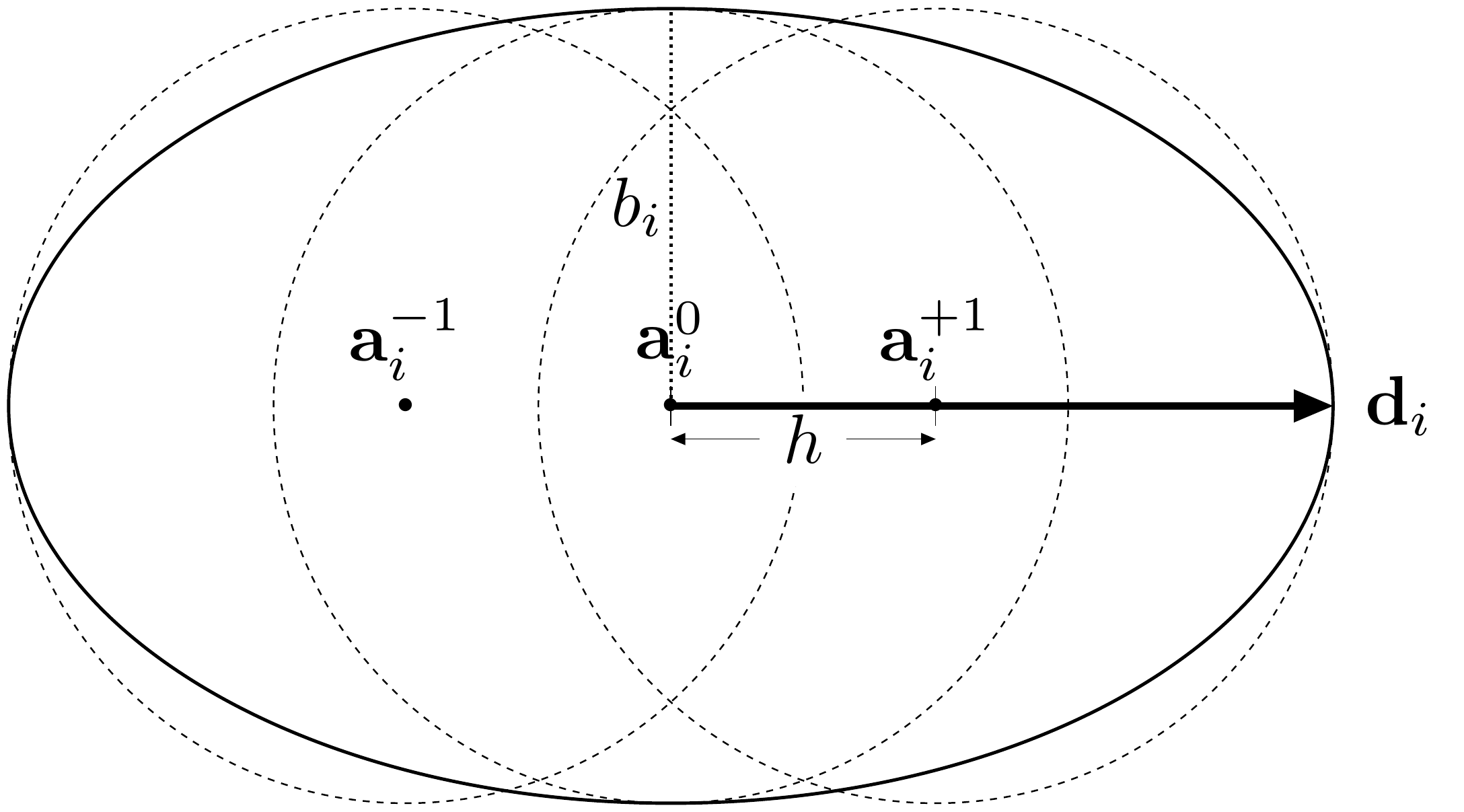}
\end{center}
  \caption{Schematic of elliptical repulsion.
   }
  \label{ellip_rep}
\end{figure}
For any elliptical particle with semi-major axis $a_i$ and semi-minor axis $b_i$, the rigid body repulsion is designed as follows. Consider three circle centers 
$\mathbf{a}_i^k = \mathbf{a}_i + kh \mathbf{d}_i,$ $h = (a_i - b_i), k\in\{-1,0,+1\}$ with radii $b_i$ (Figure~\ref{ellip_rep}) and use \cref{repul} with $c_0$ divided by 3 to calculate the forces $\mathbf{f}^k_{i}, k\in\{-1,0,+1\}$ and torque $\tau_i^{-1}$ and $\tau_i^{+1}$. 
For two elliptical particles, the total interactions are given by
\begin{equation}
\begin{gathered}
\mathbf{f}^{kl}_{ij}=\frac{c_0}{3} \frac{q}{(|\mathbf{a}_i^k - \mathbf{a}_j^l|-(b_i+b_j))^{q+1}}\frac{\mathbf{a}_i^k - \mathbf{a}_j^l}{|\mathbf{a}_i^k - \mathbf{a}_j^l|}, \quad  i\neq j,\quad k,l\in\{-1,0,+1\}, \\
\mathbf{F}^{\mathrm{rep}}_{i} = \sum_{\substack{j=1 \\  j\neq i}}^N \sum_{k,l}\mathbf{f}^{kl}_{ij}, \quad  k,l\in\{-1,0,+1\},\\
\tau_i^{\mathrm{rep}} =\sum_{\substack{j=1 \\  j\neq i}}^N \sum_{k}h (\mathbf{f}^{-1k}_{ij} - \mathbf{f}^{+1k}_{ij})\times \mathbf{d}_i,  \quad  k\in\{-1,0,+1\}.
\end{gathered}
\end{equation} 
The total repulsive potential is 
\begin{equation}
\Phi_i^{\mathrm{rep}} = \frac12\sum_{\substack{j=1 \\  j\neq i}}^N\sum_k \frac{c_0}{3} \frac{1}{(|\mathbf{a}_i^k - \mathbf{a}_j^l|-(b_i+b_j))^{q}} ,  \quad  k\in\{-1,0,+1\}.
\end{equation}
The total repulsive force is identical to \cref{repul} whenever particles are circular.


\section{Numerical Validations}
\label{sec:numerical_validations}
\subsection{Force and Torque Relations}
We validate formulas~(\ref{eq:force_torque},~\ref{eq:stress}) 
by centered difference approximation.
Following \cref{var2force}, 
\begin{equation}
\label{pos_test}
\begin{aligned}
 \mathbf{v}  \cdot \mathbf{F}_i & = -\frac{d}{d\epsilon} \Phi(\mathbf{a}_i+\epsilon \mathbf{v},\theta_i)|_{\epsilon=0} \approx -\frac{ \Phi(\mathbf{a}_i+\epsilon \mathbf{v}, \theta_i)-
  \Phi(\mathbf{a}_i-\epsilon \mathbf{v}, \theta_i)}{2\epsilon},\\
\omega \tau_i &= -\frac{d}{d\epsilon} \Phi(\mathbf{a}_i,\theta_i+\omega\epsilon )|_{\epsilon=0} \approx -\frac{ \Phi(\mathbf{a}_i, \theta_i+\omega \epsilon)-
  \Phi(\mathbf{a}_i, \theta_i-\omega \epsilon)}{2\epsilon},
  \end{aligned}
\end{equation}
where $i = 1,\dots, N$. 
We write $\Phi(\mathbf{a}_i, \theta_i)$ in place of $\Phi(\Omega,f)$ 
to emphasize that for the moment variations are taken with respect to the $i$th particle, 
while keeping the others particles fixed. 
For the three-particle setup from 
\cref{figure1}A and step size $\epsilon = 0.05,$ we get the following values:
\begin{center}
\begin{tabular}{c|c}
Centered Difference & Variations (\ref{eq:force_torque}-\ref{eq:stress}) \\
\hline
\begin{tabular}{cc}
$\mathbf{F}_i$ & $\tau_i$ \\
$\langle -0.94496,  +1.37954 \rangle $ &$+0.90685$ \\
$\langle -0.28603,  -0.46196\rangle$ & $+0.02815$\\
$\langle +1.17189, -0.90103\rangle $ &$-0.23972$
\end{tabular}
&
\begin{tabular}{cc}
$\mathbf{F}_i$ & $\tau_i$ \\
$\langle -0.83884,   +1.35038 \rangle$   &$+0.92534$ \\
$\langle -0.26879, -0.43257\rangle$     &     $+0.02923$\\
$\langle  +1.20538,  -0.91928 \rangle$ & $-0.23962$
\end{tabular}
\end{tabular}
\end{center}
The agreement between the centered difference approximation 
and the variational derivatives supports 
(\ref{eq:force_torque}-\ref{eq:stress}). 

\subsection{Single Particle Validation}

\begin{figure}
\begin{center}
\includegraphics[width=0.32\textwidth]{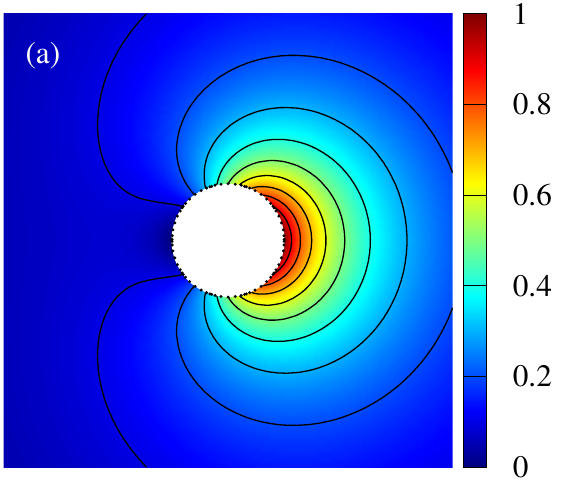}
\includegraphics[width=0.32\textwidth]{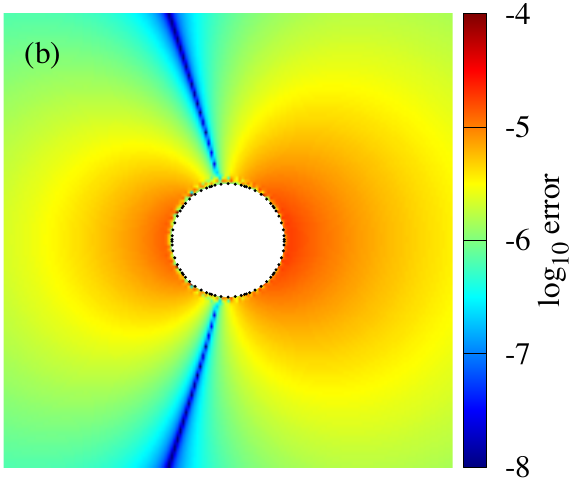}
\includegraphics[width=0.31\textwidth]{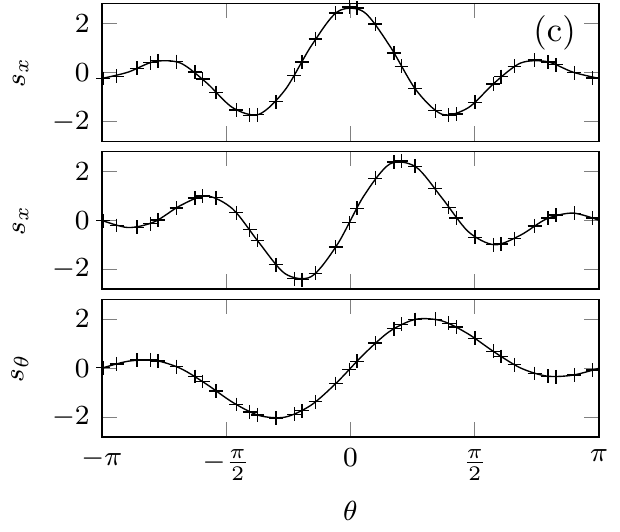}
\end{center}
  \caption{Single Janus particle simulation ($\rho=4a$): 
  (a) Surface-contour plot of the  action field solution $u.$ 
  The points along the boundary (circle) show the distribution of panels and Gauss-Legendre points.
   (b) Error plot (in logarithmic scale) for the BIE-QBX-FMM solution and analytical solution \ref{true_soln}.
  (c) From top to bottom: normal and tangential stresses calculated by BIE-QBX-FMM  ({\bf +} symbols) 
  and their corresponding analytical values (solid curves).
   }
  \label{figure5}
\end{figure}

One can analytically solve the exterior 
problem \cref{SL} for a single 
two-dimensional disk  with boundary condition \cref{JanusBC}
using the free-space Green's function \cref{Gfunc};
\begin{equation}\label{true_soln}
u(r, \theta) = \frac12\left( \frac{K_0\left(\frac{r}{\rho}\right)}{K_0\left(\frac{a}{\rho}\right)}
+\frac{K_1\left(\frac{r}{\rho}\right)}{K_1\left(\frac{a}{\rho}\right)}\cos\left(\theta\right)\right).
\end{equation}
Here $(r,\theta)$ are polar coordinates relative to the particle center (the origin),
$\mathbf{d}_1 = \langle 1, 0 \rangle$ is the particle orientation, 
and $a$ and $\rho = 4a$ are the disk radius and decay length, respectively.

\Cref{figure5}(a) shows of the action field contours for \cref{true_soln}.
Consistent with the boundary condition \cref{JanusBC},  the hydrophobic attraction is strongest in the neighborhood 
of the right semicircle. 
The smooth boundary data results in the hydrophobic attraction extending weakly to the left of the particle. 
The size of the contours 
are proportional to the decay length $\rho$,
e.g. the farthest contour $u(r,\theta)=0.1$ in \Cref{figure5}(a) would grow for larger $\rho$.  

\Cref{figure5}(b) shows the corresponding relative errors of the BIE-QBX-FMM with $N_{\rm bdy} = 70$
boundary points and QBX order $p = 6$. 
The reflectional symmetry of the error distribution is due to the symmetric particle shape and boundary condition.
The numerically computed interfacial stresses (e.g. gradients in the action field) 
are also in excellent agreement with their analytical values. 
In \Cref{figure5}(c), the + markers are for the numerically calculated pointwise 
normal and tangential stress densities, 
\begin{equation}
\label{interface_stress}
\langle s_x(\theta), s_y(\theta) \rangle  = \mathbf{T}(a,\theta) \cdot \mathbf{i}_r,\quad
s_\theta(\theta) = \mathbf{i}_r \times (\mathbf{T}(a,\theta) \cdot \mathbf{i}_r),
\end{equation}
respectively, along the particle boundary.
The smooth curves \Cref{figure5}(c) are the analytical values, obtained 
by plugging \cref{true_soln} into the integrands of the equations~\cref{eq:int}.
Thus the BIE-QBX-FMM yields a numerical solution that is highly accurate both in terms of the 
action field and its gradients along the domain boundary. 
From a physical perspective, an isolated particle has zero net force and torque (see \cref{eq:net}).
Indeed, the integrals of force and torque curves in \Cref{figure5}(c) are all zero to about eight digit accuracy.
\begin{table}[t!]
\footnotesize
\captionsetup{position=top}
\caption{Convergence tests as the QBX order $p$ and number panels per particle $N_{\rm pan}$ vary. 
The number of Gauss-Legendre points $N_{\rm GL}=6$ per panel (yielding $N_{\rm bdy}=N_{\rm GL} N_{\rm pan}$ points per particle), FMM order $p_{FMM}=10,$ 
GMRES tolerance $tol_{\rm GMRES}=10^{-13}$ and $5\times5$ computational domain  are fixed. $N_{\rm iter}$ is the number of iterations in GMRES.}\label{table1}
\begin{center}
\setlength\extrarowheight{2pt}
\begin{tabular}{ @{}c@{}@{}c@{}@{}c@{}@{}c@{}@{}c@{}| } 
\begin{tabular}{|c} 
\hline 
\\ \hline 
$N_{\rm pan}$ \\ \hline
10   \\ \hline
20 \\ \hline
40 \\ \hline
80 \\ \hline
\end{tabular}
&
\begin{tabular}{|c} 
\hline 
\\ \hline 
$N_{\rm bdy}$ \\ \hline
70   \\ \hline
140 \\ \hline
280 \\ \hline
560 \\ \hline
\end{tabular}
& 
\begin{tabular}{|@{}c@{}} \hline QBX order $p=4$\\
\hline 
\begin{tabular}{c|c}
$N_{\rm iter}$ & $l_{\infty}$ error \\ \hline
 $8$ & $3.20\times10^{-4}$ \\ \hline
  $8$ &  $2.00\times10^{-5}$ \\ \hline
  $7$ & $9.12\times10^{-7}$  \\ \hline
  $8$ &  $5.96\times10^{-8}$   \\ \hline
\end{tabular}
 \end{tabular} 
& 
\begin{tabular}{|@{}c@{}} \hline QBX order $p=6$\\
\hline 
\begin{tabular}{c|c}
$N_{\rm iter}$ & $l_{\infty}$ error \\ \hline
 $8$ &$2.01\times10^{-5}$ \\ \hline
  $7$ &  $3.93\times10^{-7}$ \\ \hline
  $8$ & $2.21\times10^{-7}$   \\ \hline
  $9$ &   $4.77\times10^{-8}$ \\ \hline
\end{tabular}
 \end{tabular} 
 & 
\begin{tabular}{|@{}c@{}|} \hline QBX order $p=8$\\
\hline 
\begin{tabular}{c|c}
$N_{\rm iter}$ & $l_{\infty}$ error \\ \hline
 $8$ & $1.23\times10^{-6}$  \\ \hline
  $7$ &  $1.23\times10^{-8}$ \\ \hline
  $8$ &  $2.84\times10^{-7}$  \\ \hline
  $8$ &  $5.01\times10^{-8}$   \\ \hline
\end{tabular}
 \end{tabular} 
\end{tabular}

\end{center}
\end{table}

Continuing with the single particle test, 
Table~\ref{table1} provides three sets of convergence tests where we tune the QBX parameters.
The purpose of these tests is to acquire a suitable parameter set  for efficient simulations.
We fix the GMRES iterative scheme tolerance  $tol_{\rm GMRES} = 10^{-13}$  
and use the FMM  to expedite the matrix-vector multiplications in GMRES iterations.
We divide each particle boundary into 
$N_{\rm pan}$ panels and fill in $N_{\rm GL}$ Gauss-Legendre points in each panel. This yields 
a total number of boundary points $N_{\rm bdy}=N_{\rm GL} N_{\rm pan}$.
In the $l_\infty$ error columns, we compute the errors with respect to the analytical solution \cref{true_soln}
over a 5 $\times$ 5 computational domain sampling at 200 $\times$ 200 cartesian grid points 
(the error excludes the values inside the particle). 
Through our setting of the underlying fast multipole order $p_{FMM}=10$, the approximation of the layer potential has about eight
digit accuracy,
leading to the observed errors `bottoming out' around that accuracy.
As a result, the results for order $p = 10$ are only marginally better than for the $p = 8$ column, but require significantly more computational time.

\begin{figure}[ht!]
\begin{center}
\includegraphics[width=0.32\textwidth]{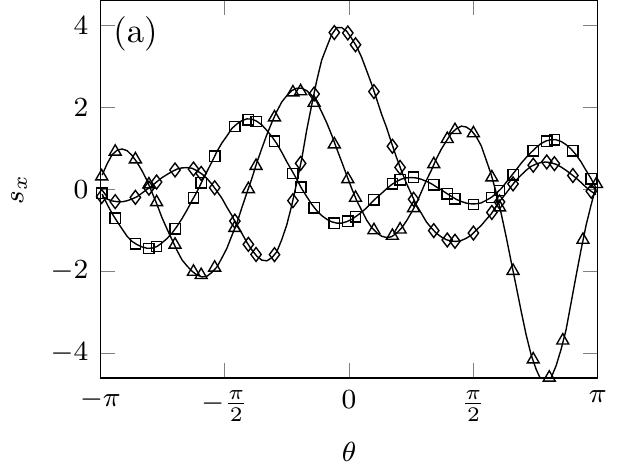}
\includegraphics[width=0.32\textwidth]{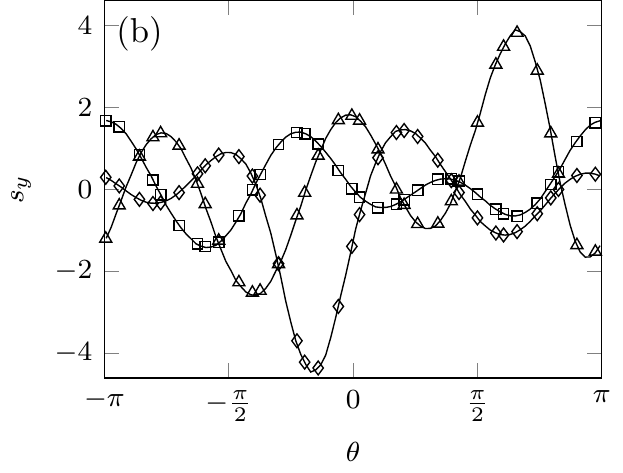}
\includegraphics[width=0.32\textwidth]{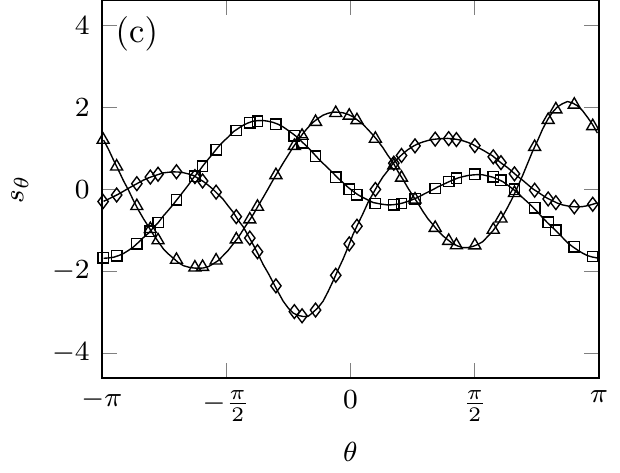}
\end{center}
\caption{Panels (a) and (b) show the normal stresses $s_x$  and $s_y,$ and panel (c) is for the tangential stress $s_{\theta}$ \ref{interface_stress}.
The symbols are: $\diamond$-diamond for particle 1, $\square$-square for particle 2 and $\triangle$-triangle for particle 3. 
Net force and torque is zero \ref{eq:net} and so the integrals of the curves in panel (a) sum to zero.
The same holds for panel (b) and for panel (c).
\label{figure6}}
\end{figure}


\subsection{Dissipative System with Constant Drag Coefficients}

As an alternative numerical scheme, the metastable final states of self-assembly particles are achievable by
using constant drag coefficients to update particle dynamics.
The updated particle dynamics of centers $\mathbf{a}_i \in \mathbb{R}^2$ and orientations $\theta_i \in \mathbb{R}$ at $t=n\Delta t$ 
are given by
\begin{equation}
{\bf a}_i^{n+1} = {\bf a}_i^{n}+\frac{1}{\xi_x}\bigg(\mathbf{F}_i + \sum_{j\neq i} \mathbf{F}_{ij}^{\mathrm{rep}}\bigg)\Delta t,\quad
\theta_i^{n+1} =\theta_i^{n} +\frac{1}{\xi_{\theta}}\tau_i \Delta t.
\label{euler}
\end{equation} 
To discretize \cref{euler}, we adopt the forward Euler scheme for configuration updates. 
We observe numerically, and use the values $\xi_{\bf x}=4\pi\mu a$ and $\xi_\theta=4\pi\mu a^3$
for an isolated circular particle of radius $a$
and $\mu=1\ $cP for water viscosity.
The numerical scheme for simulating dissipative system using proposed constant drag law is included in~\Cref{alg:qbx}.

\begin{figure}[h]
\begin{center}
\includegraphics[width=0.3\textwidth]{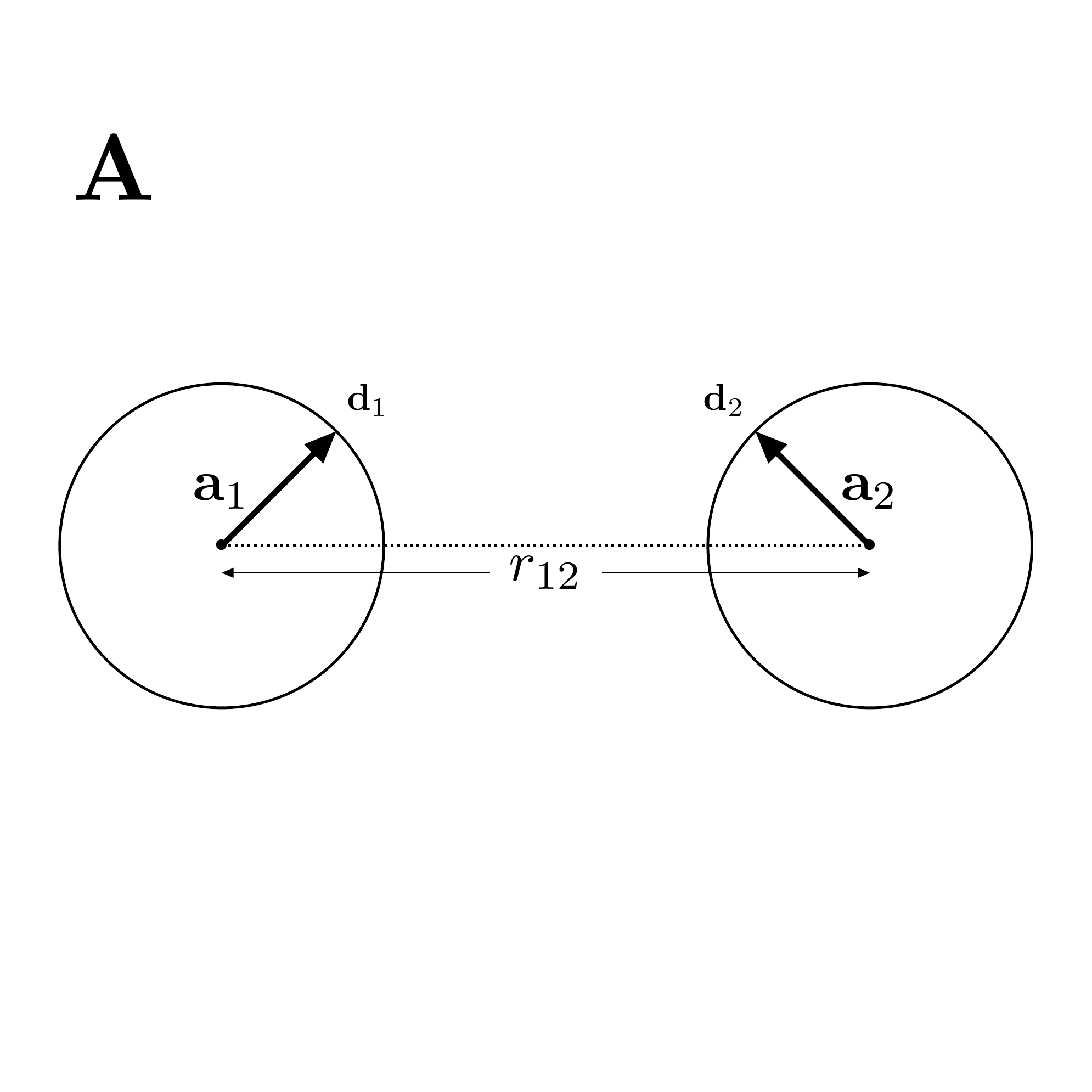}
\includegraphics[width=0.33\textwidth]{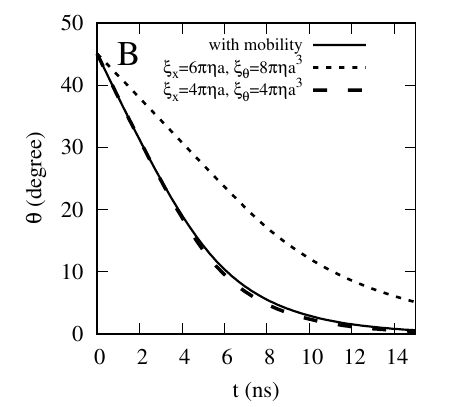}
\includegraphics[width=0.33\textwidth]{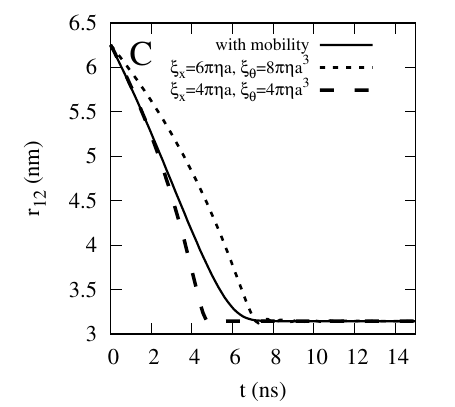}
\end{center}
  \label{N2sim}
  \caption{Panel A: Schematic of two-particle simulations; Panel B: Plot of particle 2 orientations over time in three cases; Panel C: Plot of distance $r_{12}$ over time in three cases.}
\end{figure}

This numerical test is to investigate a rough approximation of constant drag coefficients $\xi_{\bf x}$ and $\xi_\theta$. We first place two circular particles on the same horizontal axis, with centers $\mathbf{a}_1 = \langle -2.5, 0 \rangle$ and 
$\mathbf{a}_2 = \langle 2.5, 0 \rangle$, 
and orientations $\theta_1 = 45^{\circ}$ and $\theta_2 = 135^{\circ}$
(The schematic is in \Cref{N2sim}A). From the theory of HADF, the particle pair will move toward each other and rotate until the system energy reaches a minimum. Due to the effect of excluded volume repulsion, with the choice of $c_0=0.0166$ pN nm$^4$, an equilibrium distance $r_{12}$ between two particles can be measured. Three sets of simulations are performed: (1) Obtaining particle dynamics by solving a mobility problem; (2) Calculating dissipative dynamics using three dimensional translational and rotational drag coefficients $\xi_{\bf x}=6\pi\eta a$ and $\xi_{\bf x}=8\pi\eta a^3$ and (3) Calculating dissipative dynamics using translational and rotational drag coefficients $\xi_{\bf x}=4\pi\eta a$ and $\xi_{\bf x}=4\pi\eta a^3$.
Both \Cref{N2sim}B and \Cref{N2sim}C show that the dynamics obtained from case (2) have much lower 
initiative translational and rotational velocity. Case (3) gives a very good agreement in dynamics for the first few nanoseconds. To explain this finding, from Stokesian dynamics, the resistance tensor is a function of particle pair-distances and the particle resistance will be a factor of $\log(r_{ij})$ in two dimensions.
This observation shows that with a specific choice of constant drag coefficients the dynamics of many-body system may have very similar starting transition in self-assembly.

\subsection{Multiple Particle Cases}
For three-particle dynamics, 
there is no closed-form solution to compare with, 
and so we use a piecewise linear FEM to perform the numerical validations.
The three particle configuration is the same as in \ref{figure1}(a),
with centers $\mathbf{a}_1 = \langle -1, 0 \rangle,$
$\mathbf{a}_2 = \langle 1.5, 3.3 \rangle$ and 
$\mathbf{a}_3 = \langle 1.5, -1.5 \rangle,$
and orientations $\theta_1 = 18^{\circ},$ $\theta_2 = 240^{\circ}$ 
and $\theta_3 = -60^{\circ}.$ 
In the FEM solution of (\ref{JanusBC}, \ref{SL}), we truncate the unbounded exterior domain $\Omega$ 
to the box $(-20,20)\times(-20,20) \subset \mathbb{R}^2$ and apply a homogeneous Dirichlet condition on the box boundary~\cref{exp_decay}. 
To achieve an accuracy comparable to that of the BIE-QBX-FMM (Table~\ref{table1}), 
the FEM uses $N_{\rm bdy} = 250$ equally spaced points per particle boundary
and a triangular mesh with roughly $15,000$ points to discretize a truncated domain in $\Omega$.
 
\Cref{figure6}(a)--(c) compares interfacial stresses \cref{interface_stress} derived  
by BIE-QBX-FMM (empty symbols) and the FEM (solid curves).
The excellent agreement between the results suggests that
the integral equation method and the finite element method with appropriate truncation do an equally good job of calculating 
the interfacial stresses, and in practice would yield indistinguishable dynamics.  
The BIE-QBX-FMM, however, has the advantages that it uses far fewer mesh points than the FEM to achieve the same
accuracy, and that it is straightforward to discretize boundaries of moving particles using high-order quadratures. In contrast, in the FEM 
each change in particle configuration involves the generation
of a new triangular mesh as well as the artificial truncation of the domain, leading to much higher computational cost to achieve the same accuracy.

\subsection{Large Collection Simulations}
\begin{table}[ht!]
\footnotesize
\captionsetup{position=top}
\caption{Timing results for $N=\{250, 500, 1000, 2000\}$ particles with $N_{\rm bdy}=70$ 
boundary points per particle and $tol_{\rm GMRES} = 10^{-5}$.
}
\label{table2}
\begin{center}
\setlength\extrarowheight{3pt}
    \begin{tabular}{ |r| r| r| r|  r| r| r| }
    \hline
    $N_{\rm particle}$ & $N_{\rm bdy}$ &  $L_x$ & $L_y$  & Iter. & $T_{\rm GMRES}/T_{\rm total}$ & $T_{\rm total}/T_{\rm ref}$\\ \hline
        $250$ & $17500$  & $100$ & $100$  & $32$ & $0.67938$  &  $9.7$  \\
        $500$ & $35000$  &  $200$ & $200$ &  $28$ & $0.66299$  & $14.4$  \\
        $1000$ & $70000$  & $200$ & $200$  &  $32$  & $0.70062$  & $32.4$\\
        $2000$ & $140000$  &  $400$ & $400$  &  $31$  & $0.64738$  & $55.3$ \\
    \hline
    \end{tabular}
\end{center}
\end{table}
The simulation in Figure~\ref{figure3} used $N = 25$ particles and this number was sufficient
for particles to self-assemble into a vesicle shape. In realistic applications though, such as membrane fusion
or vesicle deformations, the problem is three-dimensional and the number of Janus-type  particles involved would be much larger,
on the order of thousands to tens of thousands. 
Thus we present timing results illustrating how the the evaluation of one time iteration 
\Cref{alg:qbx} scales with the particle number $N$. 

Table~\ref{table2} shows the timing results for $N = \{250,500,1000,2000\}$ particles.
The particles lie in a $L_x \times L_y$ computational domain and we use $N_{\rm bdy} = 70$
boundary points per particle.  
Their shape, disks with radius $a = 1$ and decay length $\rho = 4a$, remains 
the same as previously and their centers and orientations are randomly generated 
in a way that avoids overlapping boundaries. 

The columns include the percentage running time of GMRES (the computationally most intensive step) and 
total running time that includes the QBX initialization steps. 
In the tests of Table~\ref{table2}, which starts from random initial data, about two thirds of the simulation
time goes into solving for the surface potential $\sigma$. 
(We found that a tolerance $tol_{\rm GMRES}=10^{-5}$ gave sufficiently good numerical accuracy for the purposes of examining the particle dynamics.)
In \Cref{alg:qbx}, however, we can use the surface potential $\sigma$ 
calculated in the previous time-step as an initial guess for GMRES iterations. 
This typically reduces the GMRES iterations by a factor of four. 

The rightmost column shows the total running time relative to the reference 
time $T_{\rm ref} = 10$ sec.
for the 25 particle simulation. The results, which use an 8 core Intel(R) Xeon(R) CPU E5-2650 v4 @ 2.20GHz
for hardware, scale linearly with $N.$ 
On a modern computing cluster, most of the calculations, such as GMRES iterations, source evaluations, 
symbolic representations, FMM evaluations and numerical integrations, can run in parallel.
We therefore expect to have optimal computational cost when running large scale simulations
in future studies.

\section{Movie Captions}\mbox{} \\

\noindent
{\bf Movie S1. Three Particles} 
There are three circular particles with radius 1 centered at ${\bf a}_1=\langle0,0\rangle$, ${\bf a}_2=\langle2.5,3.3\rangle$ and ${\bf a}_3=\langle2.5,-1.5\rangle$ and the corresponding orientations $\theta_1$, $\theta_2$ and $\theta_3$ are $0.1\pi,\frac43\pi$ and $-\frac13\pi.$ In this movie, each arrow represents the director of coarse-grained lipid particles where it points from lipid head toward lipid tail. 
All white dots in the domain represent the tracers in fluid that move with respect to calculated fluid motion. The colored field from dark blue to dark red shows the magnitude of hydrophobic attraction activity and the range is from 0 to 1.
Particle 1 and 2 pair quickly aggregate and squeeze the fluid out resulted that the generated fluid flow pushes particle 3 further away from the particle pair. After few frames, due to a non-zero hydrophobic attraction activity between particles, particle 3 rotates and move toward the particle pair to reach the energy minimum. It is clear to see that
the fluid is been excluded completely at the last state of the movie.
This movies includes a total 100 time steps where the time step is $\Delta=1.0$.\\

\noindent
{\bf Movie S2. Twenty-Five Particles} 
There are 25 circular particles with radius 1 initially located on a 5-by-5 matrix grid and the initial orientations $\theta_i$ are normally distributed about $\theta=0$.
In this movie, each arrow represents the director of coarse-grained lipid particles where it points from lipid head toward lipid tail. 
All white dots in the domain represent the tracers in fluid that move with respect to calculated fluid motion. The colored field from dark blue to dark red shows the magnitude of hydrophobic attraction activity and the range is from 0 to 1.
All 25 particles begin from forming a number of micelle like groups and then assemble to three short bilayers. Here the minimal energy is not completely reached and all endpoints of 3 bilayers move toward non-zero activity field. At final equilibrium state, a vesicle is formed and a energy minimum is achieved. As suggested by HADF, the fluid is separated into two parts, outside and inside of the vesicle.
This movies includes a total 800 time steps where the time step is $\Delta t=1.0$.\\

\noindent
{\bf Movie S3. One Hundred Particles} 
This movie adopts the constant drag law to perform dissipative dynamics. We show the simulation results for 100 particle placed on a 10-by-10 grid with random orientations.
In this movie, each arrow represents the director of coarse-grained lipid particles where it points from lipid head toward lipid tail. The colored field from dark blue to dark red shows the magnitude of hydrophobic attraction activity and its range is from 0 to 1.
The parameter set is as follows, $\xi_{\bf x} = 1.5$, $\xi_\theta=2.0$ and $\Delta t = 0.5.$
All particles start from forming particle pairs or small groups then these components form 
micelles and bilayers. In order to reach energy minimum, some groups form long bilayers. 
Notice that the bilayer on the top-right corner, the transition from an arc to straight shape gives a perfect example for the process of energy minimization. Also, all micelles in the last frame have symmetric shapes.
This movies includes a total 1200 time steps.


\end{document}